\documentclass{article}

\usepackage{amsmath,amssymb,amscd,theorem,pb-diagram,boxedminipage,epic}

{\theorembodyfont{\slshape}\newtheorem{theorem}{Theorem}}
{\theorembodyfont{\slshape}\newtheorem{proposition}[theorem]{Proposition}}
{\theorembodyfont{\slshape}\newtheorem{lemma}[theorem]{Lemma}}
{\theorembodyfont{\slshape}\newtheorem{definition}[theorem]{Definition}}
{\theorembodyfont{\slshape}\newtheorem{note}[theorem]{Note}}
{\theorembodyfont{\slshape}}
{\theorembodyfont{\slshape}\newtheorem{step}{Step}}
{\theorembodyfont{\rmfamily}\newtheorem{step1}{Step}}

\newcommand{\Fq}{\ensuremath{\mathbb{F}_q}}
\newcommand{\Fp}{\ensuremath{\mathbb{F}_p}}
\renewcommand{\O}{\ensuremath{\mathcal{O}}}
\newcommand{\Ot}{\ensuremath{\widetilde{\mathcal{O}}}}
\newcommand{\G}{\ensuremath{\Gamma}}
\newcommand{\g}{\ensuremath{\gamma}}
\newcommand{\og}{\ensuremath{\bar{\gamma}}}
\newcommand{\Q}{\ensuremath{\mathbb{Q}}}
\newcommand{\Qp}{\ensuremath{\mathbb{Q}_p}}
\newcommand{\Qq}{\ensuremath{\mathbb{Q}_q}}
\newcommand{\Zp}{\ensuremath{\mathbb{Z}_p}}
\newcommand{\Zq}{\ensuremath{\mathbb{Z}_q}}
\newcommand{\R}{\ensuremath{\mathbb{R}}}
\newcommand{\C}{\ensuremath{\mathbb{C}}}
\newcommand{\Z}{\ensuremath{\mathbb{Z}}}
\newcommand{\N}{\ensuremath{\mathbb{N}}}
\newcommand{\F}{\ensuremath{\mathbb{F}}}
\renewcommand{\S}{\ensuremath{\mathbb{S}}}

\renewcommand{\k}{\ensuremath{\kappa}}

\newcommand{\ord}{\ensuremath{\text{ord}}}

\newcommand{\B}{\ensuremath{\mathcal{B}}}
\newcommand{\ep}{\hfill $\blacksquare$\vspace{\baselineskip}} 

\begin{document}

\title{Point counting in families of hyperelliptic curves}

\author{ Hendrik Hubrechts \footnote{Research Assistant of the Research Foundation - Flanders (FWO - Vlaanderen).} \\
\small{Department of mathematics, Katholieke Universiteit Leuven}\\
\small{Celestijnenlaan 200B, 3001 Leuven (Belgium)}\\
\small{\texttt{Hendrik.Hubrechts@wis.kuleuven.be}}}
\date{29 January 2007}

\maketitle

\begin{abstract}
Let $E_{\G}$ be a family of hyperelliptic curves defined by $Y^2=\bar{Q}(X,\G)$, where $\bar{Q}$ is defined over a small finite field of odd characteristic. Then with $\og$ in an extension degree $n$ field over this small field, we present a deterministic algorithm for computing the zeta function of the curve $E_{\og}$ by using Dwork deformation in rigid cohomology. The time complexity of the algorithm is $\O(n^{2.667})$ and it needs $\O(n^{2.5})$ bits of memory. A slight adaptation requires only $\O(n^2)$ space, but costs time $\Ot(n^3)$. An implementation of this last result turns out to be quite efficient for $n$ big enough.
\end{abstract}

\vspace{\baselineskip}
\section{Introduction}\label{sec:intro}
The idea that it might be interesting to compute the number of solutions to an algebraic equation over a finite field without actually trying to find these solutions is an old one, as Gauss already introduced his so called Gauss sums for it. In later times this topic led to wonderful theoretical results as for example the Weil conjectures, the pursuit of a proof of which had a very big influence on number theory and algebraic geometry. In those days computers did not yet exist as everyday objects, hence there was no real interest in computing the number of points on very concrete varieties, and for example the $\ell$-adic theory that led to the final proof of the Weil conjectures seems especially unsuited for implementations (except in the elliptic curve case).

In the eighties the idea of Koblitz \cite{KoblitzEC} and Miller \cite{MillerEC} to use elliptic curves for cryptography created a new interest in the subject, but now on this more concrete level. Later on came also the suggestion of using jacobians of hyperelliptic, $C_{ab}$ and other kinds of curves. Besides these concrete reasons the matter should of course be interesting already in itself: given a very big but finite object with an easy defining relation and lots of structure: what is its size? But there are even more reasons for investing in point counting algorithms, e.g.\ in \cite{AnshelAndGoldfeld} a one-way function is constructed which uses such an algorithm, and \cite{TsfasmanAlgebraic} describes how to use jacobians in the context of sphere packing.\\

The first interesting general algorithm that saw the light was Schoof's algorithm for calculating the number of points on an elliptic curve. This algorithm has a time complexity polynomial in the logarithm of the field size, and was optimized by Elkies and Atkin, thus resulting in the well known SEA-algorithm, see \cite{SchoofCountingPoints}. It is possible to generalize this approach to higher genus curves, but the complexity is then exponential in the genus, and this has only been done for genus equal to 2.

As higher genus curves came into view more algorithms emerged. In practice we can distinguish two kinds of problems: for a field size $p^n$ ($p$ prime), give an algorithm that works polynomially in $\log(p^n)$ for small $n$, or one that works in time $\O(np)$, hence for small $p$. The $p$-adic approach started with Satoh's canonical lift method for elliptic curves \cite{Satoh}. It was Kedlaya's famous paper \cite{KedlayaCountingPoints} that introduced the Monsky-Washnitzer cohomology in the computational world by giving a general algorithm for hyperelliptic curves in odd characteristic. These $p$-adic algorithms solve problems of the second kind, namely they are polynomial in $np$ instead of $n\log p$. But for large fields with small characteristic they have proven to be very efficient. Denef and Vercauteren generalized this algorithm to even characteristic \cite{DenefVercauteren}.

In the meantime Lauder and Wan presented an approach \cite{LauderWanCountingSmallChar} that led to an algorithm for very general varieties, and although not very practical, it works polynomially in the extension degree of the field (but exponentially in the number of variables involved). Afterwards Lauder \cite{LauderDeformation} used Dwork's deformation theory to reduce the dependency on the dimension of hypersurfaces.\footnote{A similar idea is used by Tsuzuki in the computation of Kloosterman sums \cite{TsuzukiKloosterman}.} The idea is to consider a family of hypersurfaces in one parameter $\G$, such that for example $\G=1$ gives the original problem, and $\G=0$ gives a special but easy case. Dwork's theory allows then to recover the necessary information in an efficient way by enabling to skip the most elaborate steps in Kedlaya's algorithm. We will make this more precise further on in this paper. Indeed, we develop a suggestion of Lauder to combine the Monsky-Washnitzer cohomology \`{a} la Kedlaya with a deformation. Although for curves the dependency on the dimension cannot decrease, for certain hyperelliptic curves the dependency on the extension degree $n$ will. Namely, Kedlaya's algorithm results in a time complexity\footnote{See Section 2.4 for an explanation of the notation $\Ot$.} of $\Ot(n^3)$ bit operations and $\O(n^3)$ as bit space requirements, whereas our algorithm requires respectively $\O(n^{2.667})$ and $\O(n^{2.5})$, or with a small adjustment $\Ot(n^3)$ respectively $\O(n^2)$. It is worth noting that we can find the matrix of the $p$th power Frobenius automorphism in essentially quadratic time, and only the final step, taking the characteristic polynomial of the norm of this matrix, requires the above time estimates. As we consider only odd $p$, a next project is the even characteristic case, these results can be found in \cite{HubrechtsHECEven}.

In \cite{GerkmannEC} Gerkmann has followed roughly the same kind of ideas for the elliptic curve case, including a short description for characteristic two. An implementation for Magma of this method is available on his website\footnote{\texttt{http://wwwalt.mathematik.uni-mainz.de/\~{}gerkmann/ellcurves.html}}.\\

The paper is structured as follows. We start in Section 2 with a rough sketch of the method, and a formulation of our results, Theorems \ref{thm:firstResult} and \ref{thm:secondResult}. Then we will construct the analytic setting in which the Monsky-Washnitzer cohomology with deformation lives, along with all the theoretical proofs. This theory implies the correctness of the algorithm, which is explained in Section 4, if we should work with infinite $p$-adic precision. As this is impossible, Section 5 proves that our chosen precision suffices. Section 6 computes the complexity of the algorithm and hence proves Theorems \ref{thm:firstResult} and \ref{thm:secondResult}. We then give some remarks and special interesting cases, and the final section presents results achieved with an implementation of the algorithm.\\

The author wishes to thank very much Jan Denef and Alan Lauder for proposing the topic. He also wants to thank Lauder for providing some crucial suggestions, Wouter Castryck and Joost van Hamel for the helpful discussions, Denef for the thorough reading of the paper and the correction of many small mistakes, and the anonymous referees for their many valuable suggestions and comments.
\section{Overview of the method and results}\label{sec:overview}
\subsection{Introducing notation and defining the problem} \label{sec:defining}
Let $p$ be an odd prime, $a\geq 1$ an integer and $q:=p^a$. We write $\Fq$ for the field with $q$ elements. For a field $\F$ its algebraic closure is denoted by $\bar{\F}$. Let $\Zp$ and $\Qp$ be respectively the $p$-adic integers and field, and $\Zq$ and $\Qq$ the unique unramified extensions of degree $a$ of $\Zp$ and $\Qp$. If the field $\C_p$ is a completion of $\bar{\Q}_p$, then $\C_p$ is algebraically closed as well. The corresponding valuation on these rings is written as $\ord$, normalized to $\ord\,p=1$. We will use $\sigma$ for the $p$th power Frobenius automorphism on $\C_p$. The projection $\Zq\to\Fq$ is denoted by $x\mapsto \bar x$, and the Teichm\"{u}ller lift of $\bar{x}\in\Fq$ is the unique lift $x\in\Zq$ of $\bar x$ that satisfies $x^q=x$. For a polynomial $\alpha$ in $X$ we denote with $\deg_X\alpha$ the degree with respect to $X$, and similarly we define $\deg_{\G}\alpha$. Its derivatives are written as $\alpha':=\frac{\partial}{\partial X}\alpha$ and $\dot \alpha:=\frac{\partial}{\partial\G}\alpha$.\\

Suppose we are given an equation
\[Y^2=\bar{Q}(X,\G)\qquad\text{with}\quad \bar{Q}\in\Fq[X,\G],\]
where $\bar{Q}$ is a monic polynomial of odd degree $2g+1\geq 3$ in $X$ and degree $\kappa$ in $\G$. If we suppose moreover that $\bar{Q}(X,0)$ has no double roots, then $Y^2=\bar{Q}(X,0)$ defines a hyperelliptic curve $\bar{E}_0$ of genus $g$ in Weierstrass form. Let $\og\in\bar{\F}_q$, let $n$ be such that $\Fq(\og)=\F_{q^n}$ and suppose that $\bar{Q}(X,\og)$ --- which defines the curve $\bar{E}_{\og}$ --- has no double roots. Then our goal is to compute the zeta function and hence the number of points of the curve $\bar{E}_{\og}$.

\subsection{Deformation in a nutshell} \label{sec:overDeform}
For those readers who are unfamiliar with the Monsky-Washnitzer cohomology as used by Kedlaya in order to count points on hyperelliptic curves, we refer to Kedlaya's paper \cite{KedlayaCountingPoints} or the course of Edixhoven \cite{EdixhovenCourse}\footnote{Available on
\texttt{http://www.math.leidenuniv.nl/\symbol{126}edix/oww/mathofcrypt/}.}. The idea of using deformation in the context of the Weil zeta function appears first in Dwork's paper \cite{DworkDeformation}, where he derives the resulting differential equation. The cohomology type considered in this paper is rigid cohomology, for a general definition and properties we refer to Berthelot's paper \cite{Berthelot}.

As Kedlaya we start with lifting everything to characteristic zero. Define $Q(X,\G)$ $\in\Zq[X,\G]$ as a degree preserving lift of $\bar{Q}(X,\G)$. It is clear that for some values of $\og$ in $\bar{\F}_q$ the polynomial $\bar{Q}(X,\og)$ is not squarefree, or equivalently $Y^2=\bar{Q}(X,\og)$ has affine singularities and hence does not define a hyperelliptic curve in Weierstrass form. With the resultant
\[r(\G) := \text{Res}_X\left(Q(X,\G); Q'(X,\G)\right)\]
we have that $\bar{r}(\og)\neq 0$ if and only if $\bar{Q}(X,\og)$ is squarefree, hence we will require that $\bar{r}(0)\neq 0$. We assume for a moment that the leading coefficient of $r$ is a unit in $\Zq$.
Let the ring $S$ and the $S$-module $T$ be defined by
\[S:=\Qq\left[\G,\frac 1{r(\G)}\right]^\dagger\quad\text{and}\quad
T := \Qq\left[\G,\frac 1{r(\G)},X,\frac 1{\sqrt{Q}}\right]^\dagger.\]
Note that these two structures are just variants of Kedlaya's $\Qq$ and $A^\dagger$, but with the deformation parameter $\G$ inserted in the right way. We have a differential operator $d:=\frac{\partial}{\partial X}dX:T\to TdX$, and $TdX/dT$ is a free $S$-module with basis
$\left\{\left.\frac{X^idX}{\sqrt{Q}},\frac{X^jdX}{Q} \right| i=0,\ldots,2g-1,j=0,\ldots, 2g\right\}$. Let $H_{MW}^-$ be the submodule generated by $\B := \left\{\left.\frac{X^idX}{\sqrt{Q}}\right|i=0,\ldots, 2g-1\right\}$. There exist explicit formulae ((\ref{eq:ReductionDenoms}) and (\ref{eq:ReductionXPowers}) in Section \ref{sec:freeQuotient}) that enable us to reduce elements in $H_{MW}^-$ to this basis.

On the module $H_{MW}^-$ we have the operator $F_p$, a lift of the characteristic $p$ Frobenius automorphism, and the connection
\[\nabla:H_{MW}^-\to H_{MW}^-d\G:\ f\mapsto \frac{\partial f}{\partial\G}d\G.\]
With $F(\G)$ a matrix for $F_p$ and $G(\G)$ for $\nabla$, we will prove the differential equation
\[\dot F(\G)+ F(\G)\cdot G(\G)=p\G^{p-1}G^{\sigma}(\G^p)F(\G).\]
As the zeta function of $\bar{E}_{\og}$ is completely determined by $F(\g)$, where $\g$ is the Teich\-m\"{u}l\-ler lift of $\og$, this matrix is the object that we want to compute. Indeed, the characteristic polynomial of a matrix of $F_p^{an}$, which can be derived from $F(\g)$, will be precisely the numerator of the wanted zeta function. We will see that we only need to calculate everything up to a computable finite precision in order to find this zeta function exactly.

\subsection{Overview of the algorithm}\label{sec:overAlgorithm}
We sketch briefly how the algorithm exploits the theory above. Note that we always work with $p$-adics modulo a certain power of $p$ and power series in $\G$ up to a certain power.
\begin{enumerate}\setlength{\itemsep}{0pt}
\item Lift $\bar{Q}$ to characteristic zero, and compute the resultant $r=\alpha Q+\beta Q'$.
\item Determine the matrix $G$ of the connection $\nabla$ by differentiating the basis $\B$ with respect to $\G$ and using the reduction formulae (\ref{eq:ReductionDenoms}) and (\ref{eq:ReductionXPowers}) of Section \ref{sec:freeQuotient}.
\item Deduce from $G$ the matrix $C$ whose rows form a basis of the local solutions of $\nabla=0$ by solving the differential equation $\dot{C}=-C\cdot G$.
\item Compute $\left(C^\sigma(\G^p)\right)^{-1}$.
\item Compute $F(0)$ by Kedlaya's algorithm.
\item As $F(\G)=\left(C^\sigma(\G^p)\right)^{-1}\cdot F(0) \cdot C(\G)$, we get a series expansion for $F(\G)$.
\item By representing $\F_{q^n}$ and $\Q_{q^n}$ as explained in Section \ref{sec:pAdicCalc}, compute $F(\g)$.
\item Determine the numerator of the zeta function of $\bar{E}_{\og}$ as the characteristic polynomial of
\[F(\g)^{\sigma^{an-1}}\cdot F(\g)^{\sigma^{an-2}}\cdots F(\g)^{\sigma}\cdot F(\g).\]
\end{enumerate}

\subsection{Results}\label{sec:results}
The main result of this paper is a `subcubic' time algorithm for the computation of the zeta function of certain families of hyperelliptic curves. Here we consider the extension degree $n$ to be the crucial parameter. In the formulation of these results we use the $\Ot$-notation as defined e.g.\ in \cite{ModernCompAlg}, which means that we ignore logarithmic factors. The relevant examples are that $\O(n\log n)$ and $\O(n\log n\log\log n)$ are both $\Ot(n)$. Note that we ignore the dependency on $p$, and all complexities are to be seen bitwise: for time requirements we mean bit operations, and space is also meant as number of bits. In these results the reader may take into account that a linear deformation will be most practical, which has $\kappa=1$.
\begin{theorem}\label{thm:firstResult}
Let $p$ be an odd prime and let $\bar{E}_{\og}$ be a hyperelliptic curve over $\F_{p^{an}}$ obtained by a deformation of the kind described above: $Y^2=\bar{Q}(X,\og)$ with $\bar{Q}$ defined over $\F_{p^a}$. There exists a deterministic algorithm that calculates the zeta function of $\bar{E}_{\og}$ in time
\[\Ot\left(n^{2.667}g^{6.376}a^3\k^3\right)\]
and in memory $\O\left(n^{2.5}g^5a^3\k^2\log^3g\right)$. We can also compute the same object using $\O\left(n^2g^5a^3\k^2\log^3g\right)$ memory at the cost of a time estimate $\Ot\left(n^{3}g^{6.376}a^3\k^3\right)$.
\end{theorem}
If we have done the calculation for one value of $\og$, it is faster to calculate more zeta functions from the same family.
\begin{theorem}\label{thm:secondResult}
Suppose we have as precomputation computed a sufficient approximation of the matrix of Fro\-be\-nius of a family as in Theorem \ref{thm:firstResult}. Then we can find the zeta function of a member of the family with parameter in $\F_{p^{an}}$ in deterministic time
$\Ot\left(n^{2.667}g^{5}a^3\k\right)$ and space
$\O \left(n^{2.5}g^{5}a^2\k\log g\right)$, or
$\Ot\left(n^{3}g^{5}a^3\k\right)$ respectively
$\O \left(n^{2}g^{5}a^2\k\log g\right)$.
\end{theorem}
The proof of these theorems will be given in Section \ref{sec:complexity}, after the complexity analysis.
\section{Analytic theory}\label{sec:analysis}
We have tried to make this section self-contained, to have a clear exposition of the analytic theory behind the deformation. As a consequence we will reintroduce some notation, but sometimes in a more general context, e.g.\ we will work over an algebraic closure $\bar{\F}_q$ of $\Fq$.

In \cite{KedlayaCountingPoints,KedlayaErrata}, Kedlaya presented a concrete computable form of the Monsky-Wash\-nitzer cohomology. In this section we combine this with a one-dimensional deformation as suggested by Lauder in \cite{LauderRigidCohomology}. As there are a lot of technical convergence results, we cannot present all the details of every calculation, but it will always be clear how to reconstruct those missing computations.

\subsection{Sketch of the situation}\label{sec:sketchDeform}
Let $\F_q$ be a finite field with $q=p^a$ elements, $p$ an odd prime, and suppose we are given a polynomial $\bar{Q}(X,\G)\in\F_q[X,\G]$ with $\deg_XQ=2g+1$ for some integer $g\geq 1$, and monic in $X$. For every $\og\in\bar{\F}_q$ we define the curve $\bar{E}_{\og}$ as corresponding to the equation
\[\bar{E}_{\og}\quad\longleftrightarrow\quad Y^2 = \bar{Q}(X,\og).\]
We need for our theory that $Y^2 = \bar{Q}(X,\og)$ is in Weierstrass form, hence without affine singularities. Suppose that $\bar{E}_0$ satisfies this condition, the goal is then to compute the zeta function of $\bar{E}_{\og}$ for certain `good' $\og\in\bar{\F}_q$. We remind the reader of the rings $\Zp$, $\Qp$, $\Zq$, $\Qq$ and $\C_p$ defined in Section \ref{sec:defining}. We take $Q(X,\G)\in\Zq[X,\G]$ such that $Q$ projects modulo $p$ to $\bar{Q}$, $\deg_XQ=\deg_X\bar{Q}$ and $\deg_{\G}Q=\deg_{\G}\bar{Q}$. The resultant $r(\G):=\text{Res}_X(Q(X,\G);Q'(X,\G))$ satisfies an equality
\[r(\G) = \alpha(X,\G)\cdot Q(X,\G)+\beta(X,\G)\cdot \frac{\partial Q}{\partial X}(X,\G),\quad\text{or }r=\alpha Q+\beta Q',\]
for some $\alpha(\G),\beta(\G)\in\Zq[X,\G]$. From $r(\G)$ we can find the set of good parameters
\[\S := \left\{\g\in\C_p\ \vphantom{\overline{r(\g)}}\right|\left.\ \g\text{ is a Teich\-m\"{u}l\-ler lift and }\overline{r(\g)}\neq 0\right\}.\]
The following property is obvious: for $\g\in\S$ the polynomial  $\bar{Q}(X,{\og})$ is squarefree and thus determines a hyperelliptic curve in Weierstrass form. We also have $0\in\S$ and $\ord(r(\g))=0$ for all $\g\in\S$.
\begin{definition}\label{def:constants}
Let $\rho$ be the degree of $r(\G)$, $B :=
\max\{\deg_{\G}\alpha,\deg_{\G}\beta\}$, $D :=
\max\{\deg_{X}\alpha,\deg_{X}\beta\}$ and $\kappa :=
\deg_{\G} Q $.\end{definition} It is easy to see that $\rho\leq 4g\kappa$ and that we can
choose $\alpha$ and $\beta$ such that $D\leq 2g$ and $B\leq (4g-1)\kappa$.

\subsection{The base ring S}\label{sec:ringS}
We want as base ring an extension of $\Qq$ that includes polynomials in $\G$ and ensures further on finite dimensionality of a certain quotient module. For $r(\G)=\sum_{i=0}^\rho r_i\G^i$ let $\rho'$ be the largest index for which $\ord(r_{\rho'})=0$, and define $R(\G)=\sum_{i=0}^{\rho'}r_i\G^i$. We assume for a moment that $\rho'\geq 1$, see Note \ref{not:rConstant} for the case $\rho'=0$. Define now the following overconvergent ring, equal to $\Qq[\G,\frac 1{R(\G)}]^\dagger$:
\[S := \left\{\sum_{k\in\Z}\frac{b_k(\G)}{R(\G)^k}\ \right|\left. \vphantom{\sum_{k\in\Z}\frac{b_k(\G)}{R(\G)^k}}\ b_k(\G)\in\Qq[\G],\ \deg b_k(\G)<\rho' \text{ and } \liminf_k\frac{\ord(b_k)}{|k|}>0\right\}.\]
Here the order of a polynomial is the minimum of the orders of its coefficients. It is easy to check that $1/r(\G)\in S$, and even $\sum_k b_k(\G) r(\G)^k\in S$ when $\liminf \ord(b_k(\G))/|k|>0$. We could call the $\liminf$ in the definition the `radius of convergence' of the series $s(\G)$. This is inspired by the fact that for series $\sum_i a_i\G^i\in\Qq[[\G]]$ the radius of convergence is $p^{m}$ with $m=\liminf \ord(a_i)/i$, and such a series is overconvergent if and only if $m>0$. We define the valuation of an element of $S$ as $\ord(s(\G)) := \min_k \ord(b_k)$. We can give a more exact interpretation of this overconvergence. An element $s$ of $S$ is convergent on some open disk strictly bigger than the unit disk, with finitely many smaller closed disks removed around the roots of $r(\G)$. As all these roots have norm 1, the closed disks are a subset of the unit \emph{circle}.

We could also define $S$ as the set of all rigid analytic functions $s$ on such an `overconvergent domain' $\mathcal{U}$ for which $s(\Qq\cap \mathcal{U})\subseteq \Qq$. Indeed, the completeness of $\Qq$ implies then that $s$ is defined over $\Qq$, see \cite[p.\ 196]{NonArchimedeanAnalysis}.
\begin{note}\label{not:rConstant}
The above definitions assume that $\rho'\geq 1$. If $R(\G)$ is a constant, then the theory will still work, but with a lot of changes (mostly simplifications). For example the ring $S$ will consist of all overconvergent power series $\sum_{k\geq 0}a_k\G^k$. We will not come back to this situation, as it is always clear which results require a reformulation. New proofs are nowhere necessary.
\end{note}
\begin{lemma}\label{lem:DivisionByr}
(Euclidean division for overconvergent power series) Let $f(\G)=\sum_{i=0}^\infty a_i\G^i\in\C_p[[\G]]$ and $\delta>0$, $\varepsilon\in\R$ such that $\ord(a_i)\geq \delta i+\varepsilon$ for all $i$. Then we can find $q(\G)\in\C_p[[\G]]$, $g(\G)\in\C_p[\G]$ such that $f(\G)=q(\G)R(\G)+g(\G)$ with $q(\G)=\sum_j b_j\G^j$, $\deg g<\rho'$, $\ord(g)\geq \varepsilon$ and for every $j$ we have $\ord(b_j)\geq \delta\cdot(j+\rho')+\varepsilon$.
\end{lemma}
\textsc{Proof.} It is clear that we can suppose that $R$ is monic and $\varepsilon=0$. Using Euclidean division we find for every $i$ integral polynomials $q_i(\G)$ and $g_i(\G)$ such that $a_i\G^i=q_i(\G)R(\G)+g_i(\G)$. The following properties hold: $\ord(q_i),\ord(g_i)\geq \delta i$ and $\deg q_i=i-\rho'$. The polynomials $q$ and $g$ from the lemma are now $q=\sum_i q_i$ and $g=\sum_i g_i$. For determining $b_j$ we only need $\sum_i q_i$ for those $i$ which satisfy $i-\rho'\geq j$, hence $\ord(b_j)\geq \delta(j+\rho')$.\hfill$\blacksquare$

\begin{lemma}\label{lem:SisRing}
The ring $S$ equals the set
\[\left\{\sum_{i=0}^\infty a_i\G^i+\sum_{j=1}^\infty\frac{b_j(\G)}{R(\G)^j}\ \right|\ a_i\in\Qq, \ b_j(\G)\in\Qq[\G],\ \deg b_j(\G)<\rho',\]
\[\left.\vphantom{\sum_{i=0}^\infty a_i\G^i+\sum_{j=1}^\infty\frac{b_j(\G)}{r(\G)^j}\ |\ a_i\in\Qq, \ b_j(\G)\in\Qq[\G], \deg b_j(\G)<\rho}\qquad\qquad\qquad\liminf_i\frac{\ord(a_i)}{|i|}>0\quad\text{ and }\quad\liminf_j
\frac{\ord(b_j)}{|j|}>0\right\}.\]
\end{lemma}
\textsc{Proof.} We only need to show that $f(\G):=\sum_{i=0}^\infty a_i\G^i$ with $\liminf \ord(a_i)/i>0$ can be written as $\sum_{k=0}^\infty b_k(\G)R(\G)^k$ with $\deg b_k<\rho'$ and $\liminf \ord(b_k)/k>0$, and vice versa. The latter implication is easy, for the former we use Lemma \ref{lem:DivisionByr}. We know that for all $i\geq 0$ we have $\ord(a_i)\geq \delta i+\varepsilon$ for some $\delta>0$ and $\varepsilon\in\R$. Define the sequence $f_k(\G)=\sum_i a_{i,k}\G^i$ inductively by $f_0(\G) = f(\G)$ and
\[f_k(\G) = f_{k+1}(\G)R(\G) + b_k(\G)\]
as in the previous lemma. We know then that $\ord(a_{i,k})\geq \delta i+k\delta\rho'+\varepsilon$, and hence $\ord(b_k)\geq \varepsilon+k\delta\rho'$. As $f(\G)=\sum_k b_k(\G)R(\G)^k$ we find the lemma.\ep

In order to be able to handle infinite sums of elements in $S$, we define a certain summation condition.
\begin{definition}\label{def:sigmaCondition}
Let $(s_k(\G))_{k\in\Z}$ be a sequence of elements in $S$, and $s_k(\G) = \sum_{\ell\in\Z}\frac{b_{k\ell}(\G)}{R(\G)^\ell}$. We define a set of sequences $\mathcal{S}$ as follows: $(s_k)_k$ belongs to $\mathcal{S}$ if and only if there exist $\delta>0$ and $L\in\N$ such that
\[\inf_{k\in\Z,|\ell|\geq L}\frac{\ord(b_{k\ell})}{|k|+|\ell|}>\delta.\]
\end{definition}
The following lemma is immediate.
\begin{lemma}\label{lem:StrictInequalities}
Suppose we have a sequence $(s_k)_k$ as above. Then $(s_k)_k\in\mathcal{S}$ and $\liminf_k\ord(s_k)/|k|>0$ is equivalent to the following: there exist a constant $c\in\Qq$ and $\delta>0$ such that for all $k,\ell$
\[\ord(c\cdot b_{k\ell})\geq \delta\cdot(|k|+|\ell|).\]
\end{lemma}
We now want to use the set $\mathcal{S}$ to prove convergence of an infinite sum in $S$.
\begin{lemma}\label{lem:SumSeqInS}
Let $(s_k)_{k\in\Z}$ be a sequence in $\mathcal{S}$ satisfying $\liminf \ord(s_k)/|k|>0$, then the sum $\sum_k s_k$ converges to an element of $S$.
\end{lemma}
\textsc{Proof.} The inequality in the lemma is easily seen to imply the convergence of the infinite sum, more precisely, if $s_k(\G)=\sum_\ell b_{k\ell}(\G)/R(\G)^\ell$ for every $k$, the sums $\sum_k b_{k\ell}$ converge. We suppose that the inequalities with $c$ and $\delta$ hold as in Lemma \ref{lem:StrictInequalities}, then we have for every $\ell\geq 1$ that
\[\delta\leq\inf_{k\in\Z}\frac{\ord(b_{k\ell})+\ord(c)}{|k|+|\ell|} \leq\frac{\ord(c)}{|\ell|}+\min_{k\in\Z}\frac{\ord(b_{k\ell} )}{|\ell|},\]
so that $\liminf_\ell\ord(\sum_kb_{k\ell})/|\ell|\geq \delta$.\ep

\noindent We will apply this lemma in the following technical result, where we write $R$ and $r$ for $R(\G)$ and $r(\G)$.
\begin{lemma}\label{lem:sumSatSigma}
Suppose we have for all integers $t\geq 0$ a series $s_t := \sum_{k,\ell}p_{k\ell}^{(t)}/(R^{\ell}\cdot r^k)$, a sum over $\ell\in\Z$ and $k\geq 0$, where $p_{k\ell}^{(t)}\in\Qq[\G]$, $\deg p_{k\ell}^{(t)}\leq A(k+t)$ and $\ord(p_{k\ell}^{(t)})\geq \delta(k+t+|\ell|)$ for some constants $A,\delta>0$. Then $(s_t)_t\in \mathcal{S}$ and $\sum_t s_t\in S$.
\end{lemma}
\noindent\textsc{Proof.} Fix a positive integer $N$. It is easy to verify that from $r\equiv R\bmod p$ it follows that $R^{k+N}/r^k\bmod p^N$ is an integral polynomial of degree at most $\rho N$. If we consider $s_t$ modulo $p^N$, we find that $p_{k\ell}^{(t)}\equiv 0$ as soon as $k+t+|\ell|\geq N/\delta$, so if we multiply $s_t$ with $(R\cdot r)^{N/\delta}$ we end up with a polynomial. We will now bound its degree. The worst possible degree comes from $p_{k\ell}^{(t)}$ with $k+t=N/\delta$ and $\ell=0$, which gives a degree of no more than $AN/\delta+(\rho+\rho') N/\delta$. Combining this with the above result for $R^{N+N/\delta}/r^{N/\delta}$ we conclude that $s_t\bmod p^N$ equals a polynomial of degree at most $AN/\delta+(\rho+\rho') N/\delta+\rho N=:c_1N$ divided by $R^{2N/\delta+N}=:R^{c_2N}$. If we expand the numerator as a `polynomial in $R$' we find that $s_t\bmod p^N$ can only have a nonzero coefficient for $R^m$ if $m\geq -c_2N$ and $m\leq c_1N/{\rho'}-c_2N$, or $|m|\leq \max\{c_2,c_1/{\rho'}-c_2\}N=:c_3N$. Note that the coefficient of $R^m$ also disappears if $\delta t\geq N$. If we finally write $s_t=\sum_m d_m^{(t)}/R^m$ with $\deg d_m^{(t)}<\rho'$, we find that $\ord(d_m^{(t)})\geq \max\{|m|/c_3;\delta t\}\geq |m|/(2c_3)+(\delta/2)t$. This implies that $(s_t)_t\in\mathcal{S}$ and using Lemma \ref{lem:SumSeqInS} we conclude that $\sum_t s_t\in S$.\ep

It is clear that we can always substitute a Teich\-m\"{u}l\-ler lift $\g$ from $\S$ in a series $s(\G)\in S$, because $R(\g)$ has order $0$. We conclude this section with a lemma that allows us to derive conclusions from such a substitution.
\begin{lemma}\label{lem:subsGammaInS}
Let $s(\G)=\sum_{k\in\Z}b_k(\G)/R(\G)^k\in S$ such that $\deg b_k<\rho'$ for all $k$. Suppose we have for infinitely many $\g\in\S$ that $\ord(s(\g))\geq\alpha$ for some real number $\alpha$, then also for every $k\in\Z$ we get $\ord(b_k)\geq \alpha$.
\end{lemma}
\textsc{Proof.} After multiplication by a constant we can suppose $\alpha=0$. Choose $K$ large enough to ensure that $\ord(b_k)>0$ for $|k|>K$. Define the truncated series
\[g(\G):=\sum_{k=-K}^K{b_k(\G)}{R(\G)^k},\]
so that for all the $\g\in\S$ considered we have $g(\g)\equiv f(\g)\bmod p$. We continue with contraposition. Choose $T>0$ such that $\ord(p^Tb_k)\geq 0$ for all $k$, and for at least one $k$ we have $\ord(p^Tb_k)=0$. Let $-N$ be the least of these $k$'s. Consider now
\[h(\G) := p^TR(\G)^N\sum_{k=-N}^K{b_k(\G)}{R(\G)^k}.\]
Then $\bar{h}(\G)$ is nonzero polynomial over $\F_q$, but $\bar{h}({\og})=0$ for every considered $\g\in\S$. This gives the required contradiction.\ep

We have two important consequences of this lemma. First, it will allow us to reduce to the situation described by Kedlaya if we substitute $\g\in\S$. And second, it implies that for nonzero $s\in S$ it is impossible to become zero in infinitely many $\g\in\S$.

\subsection{The ``dagger algebra'' $T$}\label{sec:moduleT}
The $\Qq$-algebra $A^\dagger\otimes\Qq$ used in Kedlaya's algorithm \cite{KedlayaCountingPoints} will be replaced by an algebra $T$ over the ring $S$. The following definition of $T$ is motivated by two requirements. First, substituting some $\g\in\S$ should reduce $T(\g)$ to the structure $A^\dagger\otimes\Qq$ considered by Kedlaya. And second, the cohomology module $TdX/dT$ should be finite dimensional. Therefore we take for $T$ the overconvergent completion of $\Qq[\G,1/R(\G),X,1/\sqrt{Q}]$ as follows:
\begin{definition}\label{def:daggerT}
\[T := \left\{ \sum_{k\in\Z}\sum_{i=0}^{2g}s_{ik}\frac{X^i}{\sqrt{Q}^k}\ \right|\left.\vphantom{\frac{X^i}{\sqrt{Q}^k}}\
\ s_{ik}\in S,\ \liminf_k\frac{\ord(s_{ik})}{|k|}>0\text{ and }(s_{ik})_k\in\mathcal{S}\right\}.\]
\end{definition}
With this definition $T$ is an $S$-algebra. An element $t(X)$ of $T$ will also be denoted by $\sum_kt_k(X,\G)/\sqrt{Q}^k$. A more exact formulation of the second aim for $T$ is as follows.
Let $\g\in\S$ and suppose ${\og}\in\F_{q^n}$ for $n$
minimal. Then $T$ with $\g$ substituted for $\G$ is precisely
$A^\dagger\otimes\Q_{q^n}$ for the curve $Y^2=\bar{Q}(X,{\og})$ over $\F_{q^n}$ as it appears in \cite{KedlayaCountingPoints}. We denote
$T(\g)$ for the structure we obtain for $T$ with the substitution $\G\leftarrow\g$, hence in fact $T(\g)=T\otimes\Q_{q^n}/(\G-\g)$.

We first prove that $A^\dagger\otimes\Q_{q^n}\subseteq T(\g)$.
Choose $\sum_{k,i}a_{ik}X^i/\sqrt{Q}^k\in
A^\dagger\otimes\Q_{q^n}$, with $a_{ik}\in\Q_{q^n}$ and
$\liminf_k\ord(a_{ik})/|k|>0$. Then $\Qq(\g)=\Q_{q^n}$ implies that we
can find for every $a_{ik}$ a polynomial $s_{ik}(\G)\in\Qq[\G]$
with the same order, degree smaller than $n$ and
$s_{ik}(\g)=a_{ik}$. This gives an element
$\sum_{k,i}s_{ik}(\G)X^i/\sqrt{Q}^k$ in $T$ which reduces to the
element we started with.

The other inclusion is easier. If $s(\G)\in S$, then $s(\g)$
converges to an element of $\Q_{q^n}$ due to the completeness of
$\Q_{q^n}$.

\subsection{A basis for the quotient module $TdX/dT$} \label{sec:freeQuotient}
A crucial r\^{o}le in the theory of Monsky and Washnitzer is
played by the co\-ho\-mo\-lo\-gy $TdX/dT$, where $d$ is the differential
operator $\frac{\partial}{\partial X}dX$ on $T$. In the
classical case the whole construction of the dagger space
$A^\dagger$ is needed to ensure this quotient has the right
properties, in particular that it is finite dimensional. Define
\[\mathcal{D}:=\left\{\frac{X^idX}{\sqrt{Q}},\ \frac{X^jdX}{Q}\ |\
i=0,\ldots, 2g-1, j=0,\ldots,2g\right\}.\] If we substitute some
$\g\in\S$ in $T$, this set $\mathcal{D}$ is a basis of $A^\dagger
dX/dA^\dagger$ as in Kedlaya's paper \cite{KedlayaCountingPoints}. Here we have a similar result.
\begin{theorem}\label{thm:dTFreeBasisB}
The $S$-module $TdX/dT$ is free with basis $\mathcal{D}$.
\end{theorem}
The proof falls down in two parts. Showing that $\mathcal{D}$ is a generating set will be proven in the following lemma, and to see that $\mathcal{D}$ is free we can work as follows. Write $\mathcal{D}=\{b_i\}$ and suppose there
exists a nontrivial relation $\sum_i s_ib_i=0$ for some $s_i\in S$,
and $s_j\neq 0$. Now Lemma \ref{lem:subsGammaInS} gives some
$\g\in\S$ for which $s_j(\g)\neq 0$, and this implies in turn the
nonzero relation $\sum_i s_i(\g)b_i(\g)=0$ in the classical case for
$A^\dagger\otimes\Q_{q^n}$ for some $n$,
contradicting the result of Kedlaya that the set $\{b_i(\g)\}_i$ is a basis for $A^\dagger\otimes\Q_{q^n}$.\hfill$\blacksquare$
\begin{lemma}\label{lem:BgeneratingSet}
$\mathcal{D}$ is a generating set.
\end{lemma}
\noindent\textsc{Proof.}
In order to express $\sum_{k\in\Z}B_k(X,\G)dX/\sqrt{Q}^k\in TdX$ as an $S$-linear combination of the elements of the basis $\mathcal{D}$, we use reduction formulae similar to those of Kedlaya in \cite{KedlayaCountingPoints}. We will first show how to reduce each $B_kdX/\sqrt{Q}^k$ to an expression in the basis $\mathcal{D}$:
\[\frac{B_kdX}{\sqrt{Q}^k}\equiv\sum_{b\in\mathcal{D}}\alpha_b^{(k)}b,\] and next we will prove that $\sum_k\alpha_b^{(k)}$ converges to an element of $S$.

In Section \ref{sec:sketchDeform} we have seen the relation $r=\alpha Q+\beta Q'$. For some numerator $B_k(X,\G)$ this becomes
\[B_k = \frac 1r(\alpha B_k)Q+\frac 1r(\beta B_k)Q'=:\frac 1r(P_kQ+R_kQ').\]
We start by reducing any expression $\sum_k B_kdX/\sqrt{Q}^k$ to a fraction with $\sqrt{Q}$ or $Q$ in the denominator. For $k\leq 0$ and $k$ even such a form is simply exact, as integrating does not change the `overconvergence' property. If $k\leq -1$ and odd we can write $\sqrt{Q}^{2\ell+1}$ as $Q^{\ell+1}/\sqrt{Q}$, and for $k> 2$ we calculate $d(R_k/{\sqrt{Q}^{k-2}})$, which gives
\begin{equation}\label{eq:ReductionDenoms}
\frac {B_kdX}{\sqrt{Q}^k}=\left(\frac 1r P_k+\frac 2{(k-2)r} R_k'\right)\frac {dX}{\sqrt{Q}^{k-2}}- 2d\left(\frac 1{(k-2)r}\frac {R_k}{\sqrt{Q}^{k-2}}\right).
\end{equation}
In both cases, where $k\leq -1$ and odd, or $k>2$, repeated application of (\ref{eq:ReductionDenoms}) gives that we end for every $B_kdX/\sqrt{Q}^k$ with a polynomial (of degree in $X$ probably very high) divided by $\sqrt{Q}$ or $Q$. The latter case gives an exact form plus a polynomial of degree at most $2g$ divided by $Q$, and for the former case we can use the following formula for $m\geq 2g$:
\[\frac{X^mdX}{\sqrt{Q}}=\frac{m-2g}{m-g+1/2} \frac{X^m-X^{m-(2g+1)}Q}{\sqrt{Q}}dX\]
\begin{equation}\label{eq:ReductionXPowers}
-\frac 1{2m-2g+1}\frac{(2g+1)X^m+X^{m-2g}Q'}{\sqrt{Q}}dX+ d\left(\frac{X^{m-2g}\sqrt{Q}}{m-g+1/2}
\right).
\end{equation}
At this point we have shown that we can always express $B_kdX/\sqrt{Q}^k$ in the basis $\mathcal{D}$, but the question remains whether the sum of these reductions converges in the right way.
We will go into details only in the most difficult case where $k$ is odd and positive. We will often use the following lemma of Kedlaya as formulated in \cite[Lemma 17.79]{CohenFrey}. The reduction in this lemma takes place in $A^\dagger$ as defined in \cite{KedlayaCountingPoints}.
\begin{lemma}\label{lem:kedlaya}(Kedlaya)
Let $h\in\Z_q[X]$ be a polynomial of degree $\leq 2g$, then for $m\in\N$ the reduction of $h(X)Y^{2m+1}dX$ (resp.\ $h(X)/Y^{2m+1}dX$) becomes integral upon multiplication by $p^\nu$ with $\nu\geq \left\lfloor\log_p\left((2g+1)(m+1)-2\right)\right\rfloor$ (resp.\ $\nu\geq \left\lfloor\log_p(2m+1)\right\rfloor$).
\end{lemma}
\begin{note}\label{not:noteAfterKedlayaLemma}
By ``the reduction'' in this lemma, the following is meant. If $$h(X)Y^{2m+1}dX=\tilde h(X)dX/Y+d(f(X,Y))$$ with $\deg h(X)\leq 2g$ and $\deg \tilde h(X)<2g$, then there exists a constant $c\in\Qq$ such that both $p^\nu \tilde h(X)$ and $p^\nu(f(X,Y)-c)$ are integral. We will always make the (allowed) assumption that $p^\nu f(X,Y)$ is integral.\\
\noindent Moreover, the proof of Kedlaya implies also the following: let $h(X)=X^k$ for some $k\leq 2g$, then with $f(X,Y)=\sum_{j\in\Z}f_j(X)/Y^{2j+1}$ and $\deg f_j(X)\leq 2g$, we  have almost always that $f_j=0$ for $j<0$ and $j>m$. The only exception is $h(X)=X^{2g}$ and $m=0$, where also $1/Y^{-1}$ can have a nonzero coefficient.
\end{note}
Because of Lemma \ref{lem:subsGammaInS} we can use Kedlaya's lemma in our theory: if it holds for every $\g\in\S$, it will also hold for $\G$. We continue with the proof of Lemma \ref{lem:BgeneratingSet}.
Clearly it suffices to consider expressions of the form $\sum_{k=1}^\infty s_kX^udX/\sqrt{Q}^{2k+1}\in TdX$, with $0\leq u\leq 2g$, $s_k=\sum_{\ell}s_{k\ell}/R^\ell$ and all $s_{k\ell}\in\Qq$. From Lemma \ref{lem:StrictInequalities} it follows that we may suppose that $\ord(s_{k\ell})\geq \delta(k+|\ell|)$ for some $\delta>0$. It is easy to see that we can multiply with an even bigger constant in order to ensure that $\ord(s_{k\ell})\geq \delta (k+|\ell|)+\lfloor \log_p(2k+1)\rfloor$ for a possibly smaller $\delta$. Writing $s_{k\ell}^{(k)}$ for $s_{k\ell}$ we use formula (\ref{eq:ReductionDenoms}) for $i=k,\ldots, 1$ as
\[\frac{s_{k\ell}^{(i)}dX}{\sqrt{Q}^{2i+1}}= \frac{s_{k\ell}^{(i-1)}dX}{\sqrt{Q}^{2i-1}}+
d\left(\frac{f_{k\ell}^{(i-1)}}{\sqrt{Q}^{2i-1}}\right),\]
where $s_{k\ell}^{(i)},f_{k\ell}^{(i)}\in S[X]$. This gives as result of the reduction process
\[\frac{s_{k\ell}dX}{\sqrt{Q}^{2k+1}}=\frac{s_{k\ell}^{(k)}dX}{\sqrt{Q}^{2k+1}}= \frac{s_{k\ell}^{(0)}dX}{\sqrt{Q}}+d\left(\sum_{i=0}^{k-1}\frac{f_{k\ell}^{(i)}}
{\sqrt{Q}^{2i+1}}\right)=: \frac{s_{k\ell}^{(0)}dX}{\sqrt{Q}}+d\varphi_{k\ell}.\]
Remember that $B :=
\max\{\deg_{\G}\alpha,\deg_{\G}\beta\}$ and $D :=
\max\{\deg_{X}\alpha,\deg_{X}\beta\}$. We will write $\deg_{\G}$ for the degree in $\G$ of the numerator in the following paragraphs.

Let us write $r^{k-i}f_{k\ell}^{(i)}= f_0+f_1Q+\cdots+f_\Delta Q^\Delta$ where $\deg_Xf_j\leq 2g$ and the $f_j$ and $\Delta$ depend on $k$, $\ell$ and $i$. As $\Delta\leq ((k-i)D+u)/(2g+1)\leq k-i+1$ we may take $\Delta=k-i+1$. We also see that \[\deg_\G f_j\leq (j+1)\kappa+\deg_\G f_{k\ell}^{(i)}\leq (j+1)\kappa+(k-i)B.\]
Write $\varphi_{k\ell}=\sum_t\varphi_{k\ell}^{(t)}/\sqrt{Q}^{2t+1}$, then the notes after Lemma \ref{lem:kedlaya} imply that we can limit ourselves to $0\leq t\leq k$. Indeed, the $\varphi_{k\ell}^{(t)}$ with $t<0$ will disappear when the reduction is finished, i.e.\ after applying formula (\ref{eq:ReductionXPowers}).
The coefficient of $1/\sqrt{Q}^{2t+1}$ in $\varphi_{k\ell}$ is
\[\sum_j \left(p_{k\ell j} := \left(f_j \text{ from the expansion of }r^{k-t-j}f_{k\ell}^{(t+j)}\right)\right)/\,r^{k-t-j},\]
where $j$ ranges from $0$ to $\lfloor(k-t+1)/2\rfloor$. Indeed, the largest $i$ for which $f_{k\ell}^{(i)}$ can appear in this sum satisfies $Q^{\Delta=k-i+1}/\sqrt{Q}^{2i+1}=1/\sqrt{Q}^{2t+1}$, hence $i=(k+t+1)/2$. We see that $\deg_{\G}p_{k\ell j}\leq (j+1)\kappa+(k-i)B\leq (2\kappa+ B)k$, as $i\leq k$ and $0\leq t\leq k$. Lemma $\ref{lem:kedlaya}$ gives that
\begin{equation}\label{eq:ordInSum}\ord\left(\sum_jp_{k\ell j}\right)\geq \ord(s_{k\ell})-\lfloor\log_p(2k+1)\rfloor\geq \delta(k+|\ell|)\geq \delta t.\end{equation}
The coefficient of $1/(R^\ell \sqrt{Q}^{2t+1})$ in $\varphi_{k\ell}$  is then $\sum_{kj}p_{k\ell j}/r^{k-t-j}$ where $k\geq t$ and $0\leq j\leq (k-t+1)/2$. In this sum we look at the terms for which $1/r^{k-t-j}=1/r^K$ for $K\geq 0$. This means that we limit ourselves to $k=K+t,\ldots,2K+t+1$, and at the same time $j$ increases from $0$ to $K=\lfloor(k-t+1)/2\rfloor$. We conclude that we have a sum $s_t = \sum_{\ell K}$ of polynomials of degree at most $(2\kappa+B)(2K+t+1)$ with valuation at least $\delta(K+t+|\ell|)$ divided by $R^\ell r^K$. Lemma \ref{lem:sumSatSigma} now implies that $(s_t)_t\in\mathcal{S}$ and hence $\sum_{k\ell}\varphi_{k\ell}\in T$.

Now we consider expressing $s_{k\ell}^{(0)}$ itself in the basis $\mathcal{D}$. For this we must reduce the $s_{k\ell}^{(0)}$ further with formula (\ref{eq:ReductionXPowers}) in order to obtain an expression of $\deg_X<2g$. As explained above, the appearing differentials are not relevant anymore. We have that $\deg_{\G} (r^k\cdot s_{k\ell}^{(0)})\leq kB$. Every time when we reduce a power of $X$ in $s_{k\ell}^{(0)}$, the degree in $\G$ increases with at most $\kappa$. If we write the reductions as $s_{k\ell}'$ this implies
\[\deg_{\G} (r^ks_{k\ell}')\leq k(B+D\kappa)+\kappa u.\]
Lemma \ref{lem:kedlaya} implies again that $\ord(s_{k\ell}')\geq \delta(k+|\ell|)$, and Lemma \ref{lem:sumSatSigma} gives the required result $s_0:=\sum_{k\ell}s_{k\ell}'\in S[X]_{<2g}$.

With a similar but easier argument we can find the same result for the reduction of $\sum_{k=1}^\infty B_kdX/Q^k$ and $\sum_{k=1}^\infty B_k\sqrt{Q}^{2k-1}dX$. We refer to \cite{HubrechtsThesis} for a complete proof.\ep

\subsection{The construction of a Frobenius lift}\label{sec:frobenius}
The purpose of this section is to construct a $p$th power Frobenius $F_p$ on $S$ and $T$. We remind the reader of the map $\sigma:\Qq\to\Qq$, the Frobenius automorphism. The definition of $F_p$ on $S$ and the module of
differentials $Sd\G$ is obvious with the following relations\footnote{It is the definition $\G\mapsto\G^p$ that forces us to take Teich\-m\"{u}l\-ler lifts for the substitution $\G\leftarrow\g$. Indeed, denote by `subs' this substitution, then we require $\text{subs}(\text{Frob}(\G))=\text{Frob}(\text{subs}(\G))$.} and
$\sigma$-linearity:
\[\G\mapsto \G^p,\qquad\frac 1{R(\G)}\mapsto\frac
1{R^\sigma(\G^p)},\qquad d\G\mapsto p\G^{p-1}d\G.\]
The fact that
$R^{\sigma}(\G^p)-R(\G)^p\equiv 0\bmod p$ gives that $1/R^\sigma(\G^p)\in S$. The definition on $T$ and $TdX$ is more complicated, namely
\[X\mapsto X^p,\qquad dX\mapsto pX^{p-1}dX,\]\[
\sqrt{Q(X,\G)}\mapsto Q(X,\G)^{\frac p2}\cdot\left(
1-\frac{Q(X,\G)^p-Q^\sigma(X^p,\G^p)}{Q(X,\G)^p}\right)^\frac
12.\] It is easy to check that this definition implies
$\left(F_p\left(\sqrt{Q(X,\G)}\right)\right)^2 =
Q^\sigma(X^p,\G^p)$. A tedious but not so hard computation shows that $F_p(\sqrt{Q(X,\G)})$ is in $T$, we omit it here. A detailed proof can be found in the author's forthcoming Ph.D.\ thesis \cite{HubrechtsThesis}.

It is clear that for a concrete $\g\in\S$ these definitions boil
down to those on $\Q_{q^n}$ and $A^\dagger\otimes\Q_{q^n}$.

\subsection{The differential equation} \label{sec:cube}
The introduction of the parameter $\G$ has as its most
important consequence the existence of the \emph{connection}
$\nabla:=\frac{\partial}{\partial \G}d\G$ on $T$. Consider the following cube of modules and morphisms:
\[\begin{diagram}
\node[3]{Td\G} \arrow[4]{e,t}{d}
\arrow[2]{s,l,-}{F_p}
  \node[4]{TdXd\G} \arrow[4]{s,l}{F_p} \\ \\
\node{T} \arrow[4]{e,t,1}{d}
\arrow[4]{s,l}{F_p}
                                  \arrow[2]{ne,l,..}{\nabla}
  \node[2]{} \arrow[2]{s}
  \node[2]{TdX} \arrow[4]{s,l,1}{F_p}
                      \arrow[2]{ne,l,..}{\nabla} \\ \\
\node[3]{T d\G} \arrow[2]{e,t,-}{d}
  \node[2]{} \arrow[2]{e}
  \node[2]{T dXd\G} \\ \\
\node{T} \arrow[4]{e,t}{d}
                                  \arrow[2]{ne,l,..}{\nabla}
  \node[4]{T dX} \arrow[2]{ne,l,..}{\nabla}
\end{diagram}  \]

It is easy to see that all faces of this cube commute. Define $H_{MW} := TdX/dT$, the first cohomology of Monsky and Washnitzer, then the above diagram implies the following commutative diagram:
\begin{equation}\label{eq:commutDiagram}\begin{CD}
H_{MW} @>{\nabla}>> H_{MW}d\G\\
@VV{F_p}V @VV{F_p}V\\
H_{MW} @>{\nabla}>> H_{MW}d\G.
\end{CD}\end{equation}
The hyperelliptic involution $\imath:T\to T:X\mapsto X,\sqrt{Q}\mapsto-\sqrt{Q}$ splits the
module $T$ in eigenspaces $T_+$ with eigenvalue 1 and $T_-$ with
eigenvalue $-1$. It is clear that $F_p\circ\imath=\imath\circ F_p$ and hence $F_p(T_-)\subseteq T_-$, and as Kedlaya points out in \cite{KedlayaCountingPoints}, we may restrict ourselves to $T_-$. The differential operators $\nabla$ and $d$ commute with $\imath$ as well, so we can restrict (\ref{eq:commutDiagram}) to $H_{MW}^- := T_-dX/d(T_-)$. From the proof of Lemma \ref{lem:BgeneratingSet} follows that $\B = \{X^idX/\sqrt{Q}\ |\ i=0,\ldots, 2g-1\}=\{b_i\}$, where $b_i:=X^idX/\sqrt{Q}$, is a basis for this free $S$-module. Write $F(\G)=(F_{i\ell})$
and $G(\G)=(G_{i\ell})$ for the matrices of
$F_p$ respectively $\nabla$ (with entries in $S$), more precisely
\[F_p(b_i)=\sum_{\ell=0}^{2g-1}F_{i\ell}b_\ell\qquad\text{ and
}\qquad \nabla(b_i)=\sum_{\ell=0}^{2g-1}G_{i\ell}b_\ell d\G.\]
Note that $F_p$ and $\nabla$ are not $S$-linear:
$F_p(\sum_is_ib_i)=\sum_iF_p(s_i)F_b(b_i)$ and
\[\nabla\left(\sum_is_ib_i\right) = \sum_i\left(\frac{\partial
s_i}{\partial\G}b_i+s_i\nabla b_i\right)d\G.\] Using the relation
$\nabla\circ F_p=F_p\circ\nabla$ on basis elements, we easily
deduce the differential equation
\[\frac{\partial}{\partial \G}F(\G)+ F(\G)\cdot G(\G) = p\G^{p-1}\cdot G^\sigma(\G^p) \cdot F(\G).\]
In order to solve this equation we can work as follows. Remember that we assumed that $\G=0$ gives a situation that can be handled by Kedlaya's algorithm, or more precisely we can compute $F(0)$ very fast as the curve $\bar{E}_0$ is defined over the small field $\Fq$. A first step is to solve the equation $\nabla=0$ locally at the origin. This means that we compute a matrix $C(\G)$ over $\Qq[[\G]]$ which satisfies $\frac{\partial}{\partial\G}C(\G)+C(\G)G(\G)=0$ and $C(0)=I$.
Indeed, with $C(\G)=(c_{i\ell})$ the vectors
$v_i:=\sum_{\ell}c_{i\ell}b_\ell$ form then a basis of the solutions of $\nabla=0$.
Next we apply $\nabla\circ F_p=F_p\circ\nabla$ to these local solutions $v_i$ in order to find that the $F_p(v_i)$ are also solutions around zero of $\nabla=0$. Hence we know that their matrix $C^{\sigma}(\G^p)\cdot F(\G)$ equals $A\cdot C(\G)$ for some constant matrix $A$. Comparing these matrices in $\G=0$ yields the equality
\[F(\G)=\left(C^\sigma(\G^p)\right)^{-1}\cdot F(0)\cdot C(\G).\]
At this point we have computed $F(\G)$ as matrix over $\Qq[[\G]]$, with entries which do not necessarily converge in Teichm\"{u}ller lifts $\g$. However, we know that $F(\G)$ is defined over $S$, hence in a last step we have to recover this representation as a matrix with entries in $S$.
\section{The algorithm}\label{sec:algorithm}
We now give a detailed exposition of the algorithm. Everything needed for a concrete implementation will be explained, either explicitly or by reference, except the handling of the $p$-adic fields and the action of Frobenius on them, which we postpone to the complexity analysis, Section \ref{sec:pAdicCalc}.

\vspace{\baselineskip}
\noindent\textsc{Input.} A prime number $p\geq 3$, the field $\Fq$ with $q=p^a$ and represented as $\Fp[x]/\bar{\chi}(x)$, $\og$ and the field $\Fq(\og)$ represented as $\F_{q^n}=\Fq[y]/\bar{\psi}(y)$ of extension degree $n$ over $\Fq$. A polynomial $\bar{Q}(X,\G)$ over $\Fq$, monic in $X$ of degree $2g+1$ where $g\geq 1$, such that $\bar{Q}(X,0)$ and $\bar{Q}(X,\og)$ are both squarefree.

\noindent\textsc{Output.} The zeta function of the smooth completion of the curve $Y^2=\bar{Q}(X,\og)$.

\begin{step}\label{step:1}
Some preliminary computations.
\end{step}
Let $\chi(x)$ be a naive lift to characteristic zero of $\bar{\chi}(x)$, such that the coefficients of $\chi(x)$ are small integers, e.g.\ from the set $\{-\frac{p-1}2,\ldots,\frac{p-1}2\}$, and so that $\deg\chi=\deg\bar{\chi}$. We represent $\Qq$ as $\Qp[x]/\chi(x)$. Lift $\bar{Q}$ to $Q$ by lifting its coefficients in the same way, and set $\k:=\deg_{\G}Q$. Define now the following constants: $\eta := \lceil 2g\log_pg+g\rceil$,
\begin{alignat*}{2}
N_0 &:= \left\lceil nga/2+(2g+1)\log_p2\right\rceil,&\quad
N_8 &:= an\left\lfloor\log_p(g)+2\right\rfloor+\lfloor 2gan(\log_pg+3)\rfloor,\\
N_b &:= N_0+N_8,&\quad N_{\G} &:=(2N_b+5)(8g+2)\kappa p+1,\\
M &:= p(2N_b+4)+(p-1)/2, &\quad N_3 &:= (2\eta +1)\lceil \log_pN_{\G}\rceil,\end{alignat*}
\[N_4 := \left\lceil\log_2(N_{\G})\eta\left(2\left\lceil\log_p (N_{\G}/p)\right\rceil +\lceil\log_pN_{\G}\rceil\right)\right\rceil,\]
and finally we have
\begin{alignat*}{2}
N_6 &:= 3\eta\lceil\log_pN_{\G}\rceil, &\qquad N_a &:= N_b + N_3+N_4+N_6.
\end{alignat*}

We work always modulo $\G^{N_{\G}}$. Moreover, in steps 2 to 6 we work modulo $p^{N_a}$ and in the last two steps modulo $p^{N_b}$. $N_0$ is the precision needed in order to recover the zeta function correctly, and $N_i$ for $i\in\{3,4,6,7\}$ is the precision lost in step $i$. Compute $\alpha(X,\G)$, $\beta(X,\G)$ and the resultant $r(\G)$ such that $r(\G)=\alpha(X,\G)\cdot Q(X,\G)+\beta(X,\G)\cdot Q'(X,\G)$. This is a matter of simple linear algebra, see e.g.\ Section 6.3 in \cite{ModernAlgebra}. We note that we need the representation of $\alpha$ and $\beta$ in their classical form in order to obtain the lowest degree in $X$. Define finally $C:=\max\{\deg_{\G} \alpha,\deg_{\G}\beta\}$ and $D:=\max\{\deg_{X} \alpha,\deg_{X}\beta\}$.

\begin{step}\label{step:2}
The matrix $G(\G)=H(\G)/r(\G)$ of the connection.
\end{step}
We compute the matrix $H(\G)$ as explained in the proof of Proposition \ref{thm:boundsOnH}.

\begin{step}\label{step:3}
The local solutions $C$.
\end{step}
The matrix $C(\G)$ satisfies the condition $C(0)=I$ and $\dot{C}+C\cdot G=0$, which becomes $r\cdot \dot{C}=-C\cdot H$. Write $H=\sum_{i=0}^hH_i\G^i$ and $C=\sum_{i=0}^{N_{\G}-1}C_i\G^i$ where the $H_i$ and $C_i$ are matrices over $\Qq$ and $C_0=I$. With $r(\G)=r_\rho\G^\rho+\cdots+r_1\G+r_0$ we obtain
\[\sum_{k=1}^{N_{\G}-1}\left(\sum_{i=0}^{\min\{\rho,k\}}(k-i)r_iC_{k-i}\right) \G^{k-1} = \sum_{k=0}^{N_{\G}-2}\left(\sum_{i=0}^{\min\{h,k\}}-C_{k-i}\cdot H_i\right)\G^k.\]
Proceeding recursively, the resulting formulae are (this proves at once the uniqueness of $C$)
\begin{align*}
C_1&=-\frac 1{r_0}C_0H_0,\text{\quad similar $C_2$, \ldots;\quad for $k\geq \rho$ and $k\geq h$ we have}\\
C_k &= \frac 1{kr_0}\left(-C_{k-1}H_0-\cdots-C_{k-h-1}H_h-r_1(k-1)C_{k-1}-\cdots\right.\\
&\hphantom{=}\qquad\qquad\qquad\left.-r_{\rho}(k-\rho)C_{k-\rho}\right).
\end{align*}

\begin{step}\label{step:4}
$\left(C^\sigma(\G^p)\right)^{-1}.$
\end{step}
First compute $C^\sigma(\G^p)$, then invert the resulting matrix using quadratically convergent Newton approximation as in \cite[Section 5.2.2]{LauderDeformation}.

\begin{step}\label{step:5}
Frobenius for the $\G=0$ case.
\end{step}
Compute $F(0)$ with Kedlaya's algorithm, but with the (higher) $p$-adic precision $N_a$.

\begin{step}\label{step:6}
Frobenius for the general case.
\end{step}
Calculate first $F(\G)= \left(C^\sigma(\G^p)\right)^{-1}\cdot F(0)\cdot C(\G)$ modulo $\G^{N_{\G}}$, and then $r(\G)^MF(\G)$ modulo $\G^{N_{\G}}$.

\noindent We now switch to the precision $p^{N_b}$.
\begin{step}\label{step:7}
Frobenius for the concrete $\og$.
\end{step}
Compute (e.g.\ as in \cite{ShoupMinimalPolynomial}) the minimal polynomial $\bar{\varphi}(y)$ of $\og$ over $\Fq$, and write $\F_{q^n}$ as $\Fq[y]/\bar{\varphi}(y)$. Lift $\bar{\varphi}(y)$ to $\varphi(y)$ over $\Qq$ such that $\Q_{q^n}=\Qq[y]/\varphi(y)$ and $\g=y$ is a Teichm\"{u}ller lift as explained in Section \ref{sec:pAdicCalc}. This trick\footnote{This trick was suggested to us by Alan G.B.\ Lauder, personal communication.} allows us to evaluate an entry $g(\G)/r(\G)^M$ of $F(\G)$ by simply reducing $g(\g)=g(y)$ modulo $\varphi(y)$, and then multiplying with $(1/r(y))^M$ modulo $\varphi(y)$. All these calculations can be done very quickly as explained in \cite{BernsteinFastMultiplication}.

\begin{step}\label{step:8}
The zeta function.
\end{step}
The matrix $\mathcal{F}$ of $F_p^{an}$ equals
\[\mathcal{F}:=F(\g)^{\sigma^{an-1}}\cdot F(\g)^{\sigma^{an-2}}\cdots F(\g)^{\sigma}\cdot F(\g).\]
This product can be calculated as in \cite{KedlayaCountingPoints}: $M_1:=F(\g)^\sigma\cdot F(\g)$, $M_2 := M_1^{\sigma^2}\cdot M_1$, $M_3 := M_2^{\sigma^4}\cdot M_2$ etc, and taking the product of those $M_i$ implied by the binary expansion of $an$. As shown by Kedlaya in \cite{KedlayaCountingPoints}, the zeta function $Z(t)$ of (the smooth completion of) the hyperelliptic curve $\bar{E}_{\og}:Y^2=\bar{Q}(X,\og)$ is given by
$Z(t)=\det\left(I-\mathcal{F}t\right)(1-t)^{-1}(1-q^{an}t)^{-1}$, and using Newton's formula
\[\det(I-\mathcal{F}t)=\exp\left(-\sum_{k=1}^\infty\text{Tr}(\mathcal{F}^k)\frac{T^k} k\right),\]
we compute $\det(I-\mathcal{F}t)$ modulo $p^{N_0}$. As Kedlaya showed in \cite{KedlayaCountingPoints}, each coefficient of $a_i$ of this polynomial satisfies $|a_i|\leq 2^{2g}p^{ang/2}$, which allows us finally to recover the zeta function.
\section{Proof of correctness}\label{sec:correct}
\subsection{The behavior of the matrix of Frobenius}\label{sec:estimatesF}
In order to determine the required accuracy in our algorithm, we need precise estimates on the convergence rate of the entries of the matrix $F(\G)$. This will be investigated in the following proposition. Define $\mu:=(pg-4)/2$.
\begin{proposition}\label{thm:estimatesF}
Let $N\in\N$ and $f(\G)$ an entry of $F(\G)$, reduced modulo $p^N$. If $N\geq \mu$, then with $\chi_1 = p(2N+4)+(p-1)/2$ and $\chi_2=(2N+5)(8g+2)\kappa p+1$ we have that $r(\G)^{\chi_1}f(\G)$ is a polynomial of degree at most $\chi_2$. Moreover, $\ord(f(\G))\geq-(\log_p g+2)$.
\end{proposition}
\textsc{Proof.} The matrix $F(\G)$ is obtained by computing Frobenius on a basis element $X^k/\sqrt{Q}$ and then reducing this result using formulae (\ref{eq:ReductionDenoms}) and (\ref{eq:ReductionXPowers}):
\[F_p\left(\frac{X^k}{\sqrt{Q}}\right)=\sum_{i\geq 0}{-\frac 12\choose i}\frac{X^{kp}\left(Q^\sigma(X^p,\G^p)-Q(X,\G)^p\right)^i}{Q(X,\G)^{pi+\frac p2}}=:\sum_{i\geq 0}\frac{B_i}{Q^{pi+\frac p2}}.\]
At this point we know that $\ord(B_i)\geq i$, $\deg_X B_i \leq kp+i((2g+1)p-1)$ and $\deg_\G B_i\leq \kappa pi$. Let $\alpha_i := pi+\frac{p-1}2$. In a first stage we use formula (\ref{eq:ReductionDenoms}) as often as needed in order to find $B_i'/\sqrt{Q}$ as reduction of $B_i/Q^{\alpha_i}$:
\[F_p\left(\frac{X^k}{\sqrt{Q}}\right)dX\equiv\sum_{i\geq 0}\frac{B_i}{\sqrt Q}dX.\]
Formula (\ref{eq:ReductionDenoms}) has to be used $\alpha_i$ times for $B_i$, hence we find \[\deg_XB_i'\leq\deg_XB_i+\alpha_iD;\qquad \deg_{\G}(r^{\alpha_i}B_i')\leq \deg_{\G}B_i+\alpha_i B.\]
The next step consists of applying formula (\ref{eq:ReductionXPowers}) in order to decrease the degree in $X$. If we write the result as $B_i'/\sqrt{Q}dX\equiv B_i''/\sqrt{Q}dX$ with $\deg_XB_i''<2g$, then
\[\deg_{\G}(r^{\alpha_i}B_i'')\leq\deg_{\G}(r^{\alpha_i}B_i')+ (\deg_X(r^{\alpha_i}B_i') -(2g-1))\kappa.\]
At this point we only need to know which $B_i''$ will not be zero modulo $p^N$. Using Euclidean division we can write $B_i(X,\G)=t_i(X,\G)Q(X,\G)^{\alpha_i+1}+s_i(X,\G)$, with $\deg_X s(X,\G)<(\alpha_i+1)(2g+1)$. This implies that
\begin{equation}\label{eq:splitBi}\frac{B_i}{Q^{\alpha_i+1/2}}= t_i\sqrt{Q}+\frac{s_i}{Q^{\alpha_i+1/2}},\end{equation}
where writing the second term as $\sum_j s_{ij}/\sqrt{Q}^{2j+1}$ would give $s_{ij}=0$ for $j<0$. It is easily seen that for such an expression Lemma \ref{lem:kedlaya} remains true, regardless of the condition $\deg h(X)\leq 2g$, and hence for the valuation of the reduction of $s_i/Q^{\alpha_i+1/2}$ we find at least
\begin{equation}\label{eq:estimateValuation} i-\lfloor\log_p(2\alpha_i+1)\rfloor\geq i-(1+\log_p(2i+1))\geq i/2-2.\end{equation}
This means we can confine ourselves to those $i$ for which $i/2-2\leq N$, or $i\leq 2N+4$. The resulting denominator $r^{\alpha_i}$ gives then $\chi_1 = \alpha_{2N+4}$, and the degree in $\G$ of the numerator satisfies
\[\deg_{\G}r^{\alpha_{2N+4}}B_{2N+4}''\leq (2N+4)(2\kappa p + 8g\kappa p)+5g\kappa p=:\chi_2-1,\]
by a tedious but trivial calculation. We do have however not yet taken the reduction of $t_i\sqrt{Q}$ in (\ref{eq:splitBi}) into account. We will show that the condition $N\geq \mu$ ensures that this contribution can be ignored. It is immediately clear that as $k\leq 2g-1$,
\begin{equation}\label{eq:inequality1}\deg_Xt_i\leq \deg_XB_i-(2g+1)(\alpha_i+1)\leq pg-i.\end{equation}
If $t_i\neq 0$ then $i\leq pg$ by (\ref{eq:inequality1}), and as we required that $2N+4\geq pg$ in the theorem, these $t_i$ will not be responsible for a bigger degree in $\G$ than $\chi_2-1$.

In order to show that $\ord(f(\G))\geq -(\log_p(g)+2)$ we have to do a little more work. The inequality certainly follows from (\ref{eq:estimateValuation}) for the part $s_i/Q^{\alpha_i+1/2}$ in (\ref{eq:splitBi}) above, but for $t_i\sqrt{Q}$ this is not so clear. Writing $t_i=t_{i0}+t_{i1}Q+\cdots+t_{i\Delta_i}Q^{\Delta_i}$ we find that $\Delta_i\leq \frac{pg-i}{2g+1}$. Lemma \ref{lem:kedlaya} implies then that the valuation of the reduction of $t_i\sqrt{Q}$ is dominated by the possible value for $t_{i\Delta_i}Q^{\Delta_i}$, being
\[i-\lfloor\log_p\left((2g+1)(\Delta_i/2+1)\right)-2\rfloor\geq -(\log_pg+2),\]
as can be checked directly.\ep

\subsection{Estimates on $H$ and $C$}\label{sec:estimatesHC}
\begin{proposition}\label{thm:boundsOnH}
The matrix $H=rG$ becomes integral after multiplying by $p^{\frac{10g}{p-1}}$ and consists of polynomials in $\G$ of degree at most $8g\kappa$.
\end{proposition}
\textsc{Proof.} The entries $(rG_{i\ell})$ of $H$ are obtained by computing $$r\cdot\nabla\left(\frac{X^idX}{\sqrt{Q}}\right) = -\frac r2\frac{X^i\dot QdXd\G}{\sqrt{Q}^3}\equiv -\frac 12\frac{(P+2R')}{\sqrt{Q}}dXd\G,$$where $P:=X^i\alpha\dot{Q}$ and $R:=X^i\beta\dot{Q}$ and we have used formula (\ref{eq:ReductionDenoms}). This has to be reduced further with formula (\ref{eq:ReductionXPowers}). We have that
\begin{align*}
\deg_X(P+2R')&\leq 2g-1+D+2g\leq 6g-1,\text{ and}\\
\deg_\G(P+2R')&\leq B+\kappa-1\leq 4g\kappa.
\end{align*}
Each reduction of $X^mdX/\sqrt{Q}$ using (\ref{eq:ReductionXPowers}) increases $\deg_\G$ by at most $\kappa$, decreases $\deg_X$ by at least 1, and introduces a denominator $2m-2g+1$. Together this gives a denominator $\prod_{m=6g-1}^{2g}(2m-2g+1)$, easily seen to be a divisor of $(10g)!$, with order at most $\frac {10g}{p-1}$. The degree in $\G$ is obtained by adding at most $4g-1$ times $\kappa$ to $4g\kappa$.\ep

\noindent We note that if we would use Lemma \ref{lem:kedlaya}, the above estimate could be improved, but this would not have any influence on the complexity estimates in the end.

In order to control the valuation of the entries of $C$, we need bounds for $F(\G)$ and $F(0)^{-1}$. The former was obtained in Proposition \ref{thm:estimatesF}, and for the latter we have the following lemma.

\begin{lemma}\label{lem:ordCInvers} The matrix $F(\G)^{-1}$ is also defined over $S$, and if $\ord(F(\G))\geq \varepsilon$ for some real number $\varepsilon$, then $\ord(F(\G)^{-1})\geq (2g-1)\varepsilon-g$.
\end{lemma}
\textsc{Proof.} Define $d(\G):=\det(F(\G))$, then clearly $d(\G)\in S$. Let $\g\in\S$ and choose $m$ such that $\og\in\F_{p^m}$, then the Weil conjectures imply that
\[\prod_{i=0}^{m-1}\sigma^i(d(\g))=p^{gm},\]
which gives immediately that $\ord(d(\g))=g$. Lemma \ref{lem:subsGammaInS} gives then the same valuation for $d(\G)$. From the following lemma we can see that $d(\G)^{-1}\in S$, and it is clear that $\ord(d(\G)^{-1})=-g$. From linear algebra we know that an entry of $\alpha(\G)^{-1}$ is a minor of $\alpha(\G)$ multiplied by $d(\G)^{-1}$, which concludes the proof.\hfill$\blacksquare$
\begin{lemma}\label{lem:inverseFunction}
Let $s\in S$ and $\varepsilon\in\R$ such that for all $\g\in\S$ we have $\ord(s(\g))=\varepsilon$. Then $1/s\in S$.
\end{lemma}
\textsc{Proof.} Multiplying $s$ with a constant if necessary we may suppose that $\varepsilon = 0$, and hence $\ord(s(\G))= 0$ as well. Let $K$ be a finite extension field of $\Qq$ such that $R(\G)$ splits completely in $K$. If we can prove that $1/s\in S\otimes K$, then trivially also $1/s\in S$.

We write $s(\G)=\sum_{k\in\Z}b_k(\G)R(\G)^k$ and choose $\delta>0$ such that $\ord(b_k(\G))>\delta\cdot |k|$ for $k\gg 0$. This implies that $b_k(\G)\equiv 0\bmod p^{\delta N}$ for $|k|\geq N$ and $N\gg 0$. From now on we only consider $N$ big enough for these inequalities to hold. Let $s_N(\G):=(s(\G)\bmod p^{\delta N})$, then it is clear that $s_N^{(1)}(\G):=R(\G)^Ns_N(\G)$ is a polynomial of degree less than $2\rho'N$. Let $\mathcal{R}$ be the set of roots of $R(\G)$ in $K$. As $\ord(s_N^{(1)}(\g))=0$ for all $\g\in\S$, there exists a factorization (with $c\in\Z_q^\times$)
\[s_N^{(1)}(\G)\equiv c\prod_{x\in\mathcal{R}}(\G-x)^{\ell_x}\bmod p,\quad\sum_{x\in\mathcal{R}}\ell_x< 2\rho'N.\]
If we define the integral polynomial $s_N^{(2)}(\G) := R(\G)^{2\rho'N}/\prod (\G-x)^{\ell_x}$, then we can also find a rational number $t>0$ and an integral polynomial $s_N^{(3)}(\G)$ such that
\[s_N^{(1)}(\G)=c\frac{R(\G)^{2\rho'N}}{s_N^{(2)}(\G)}-p^ts_N^{(3)}(\G).\]
Inverting $s_N(\G)$ and substituting the above expressions yields
\[\frac 1{s_N(\G)}=\frac{R(\G)^N}{s_N^{(1)}(\G)}=c^{-1}R(\G)^{N-2\rho'N}s_N^{(2)}(\G) \left(\sum_{k=0}^\infty\left(p^t\frac {s_N^{(2)}(\G)s_N^{(3)}(\G)} {cR(\G)^{2\rho'N}}\right)^k\right).\]
Now working as at the end of the proof of Lemma \ref{lem:sumSatSigma} we see that $(1/s_N)_N$ converges to an element of $S$.\ep

In step 3 of the algorithm we compute $C(\G)=\sum_k C_k\G^k$ using an expression of the form $C_k=\frac 1k(\cdots)$. With a trivial argument this yields an $\O(k)$-bound for $-\ord(C_k)$. It is however possible to do much better if we use Dwork's trick. This works as follows: suppose $C(\G)$ converges on a certain disk with radius $\varepsilon<1$. Then as $C(\G)=F(0)^{-1}C^\sigma(\G^p)F(\G)$, $C$ will also converge in the disk with radius $\sqrt[p]{\varepsilon}$. Repeating this process leads to convergence on the open unit disk, and if we keep track of the orders we find even an explicit bound, as proven in the following proposition.

\begin{proposition}\label{thm:boundsOnC}
Write $C(\G)=\sum_{k=0}^\infty C_k\G^k$, then we have
\[\ord(C_k)\geq-\lceil \log_p(k+1)\rceil (2g\log_p(g)+g).\]
\end{proposition}
\textsc{Proof.} As just mentioned we use the deformation relation for $F(\G)$. Write $\alpha$, $\beta$ for lower bounds for the order of $F(\G)$, $F(0)^{-1}$. If we look at the coefficients of $\G^0$ up to $\G^{p-1}$ (as matrices), we see that
\[\ord(C_k)\geq \ord(C_0^\sigma)+\alpha+\beta\]
for $k=0,\ldots, p-1$, where $\ord(C_0^\sigma)=\ord(C_0)=0$. Repeating this gives
\[\ord(C_k)\geq\max_{j=0}^{p-1}\ord(C_j)+\alpha+\beta\]
for $k=p,\ldots, p^2-1$, so e.g.\ $\ord(C_{p^2-1})\geq 2(\alpha+\beta)$. Pushing this further gives in general
\[\ord(C_k)\geq\lceil\log_p(k+1)\rceil(\alpha+\beta).\]
We have seen that we can take $\alpha=-(\log_p(g)+2)$ and $\beta=-(2g-1)(\log_p(g)+2)-g$, which proves the proposition.\ep

\noindent By a similar calculation we find with $(C^\sigma(\G^p))^{-1}=\sum_{k\geq 0}D_k\G^k$ that
\[\ord(D_k)\geq -\lceil \log_p(k/p+1)\rceil (2g\log_pg+g).\]

\subsection{Proof of correctness}\label{sec:subsecCorrect}
In this section we will prove that the chosen values for the $N_i$ and $M$ in the algorithm are sufficient in order to find the correct result. First of all, just as in Kedlaya's paper \cite{KedlayaCountingPoints} it suffices to compute the zeta function modulo $p^{N_0}$ in order to determine it exactly. During the algorithm every multiplication of non integral $p$-adic numbers can generate a decrease in the accuracy. We will handle this problem by using Lemma 22 from \cite{LauderDeformation} that says that the introduced error in the multiplication of two matrices with negative $p$-orders $-x$ and $-y$ is at most $x+y$. We will show that the loss in accuracy in step $i$ is $N_i$ for $i\in\{3,4,6,8\}$.

The fact that we have to compute $F(\g)$ with precision $N_b = N_0+N_8$ can be seen as follows. The error that could be introduced while computing the product of the matrices $F(\g)^{\sigma^i}$ has valuation at most $an(\log_pg+2)$, and computing $\det(I-\mathcal{F}t)$ cannot remove more than $\frac{2g}{p-1}+\log_p(2g)+2gan(\log_pg+2)\leq 2gan(\log_pg+3)$ of the resulting accuracy. Indeed, the worst appearing denominators in Newton's formula are $(2g)!$ and $2g$, and $\ord(\text{Tr}(\mathcal{F}^{2g}))\geq -2gan(\log_pg+2)$. Proposition \ref{thm:estimatesF} shows that the chosen $N_{\G}$ and $M$ suffice for step 7, and Lemma \ref{lem:boundsErrorC} will give $N_3$. In the formula from \cite{LauderDeformation} for $(C^{\sigma}(\G^p))^{-1}$ we have to calculate products of the form $D_kCD_k$ where $D_k$ is the $k$th approximation of $(C^{\sigma}(\G^p))^{-1}$, obtained during the $k$th iteration of the Newton approximation. This has to be done $\log_2N_{\G}$ times, hence we get $N_4$. We can use the estimates for $C$, $F(0)$ and $(C^{\sigma}(\G^p))^{-1}$ for step \ref{step:6}, which is dominated by $N_6$.

\noindent A naive estimate using the formulae in step 3 would give a huge bound on the loss in precision, but the following lemma shows that we can do better.
\begin{lemma}\label{lem:boundsErrorC}
Denoting the exact solution of $\nabla=0$ by $C$ and our computed solution modulo $p^{N_a}$ by $\tilde C$ we find for $E=\sum_{k\geq 0}E_k\G^k:=p^{-N_a}(C-\tilde C)$ that
\[\ord(E_k)\geq-(4g\log_p(g)+2g+1)\lceil\log_p(k+1)\rceil.\]\end{lemma}
\textsc{Proof.} This lemma is from \cite{LauderRecursive}, but we give a slightly different proof.\\
Truncating the right hand side in the computation of $\tilde C_i$ modulo $p^{N_a}$ in step 3 implies that $\tilde C$ is a solution of $r\dot{\tilde C}+\tilde CH=p^{N_a}\mathcal{E}$ with $\mathcal{E}$ an integral matrix. We choose a matrix $K$ such that $E+KC=C$. Then as $\nabla C=0$ and $H=rG=-rC^{-1}\dot{C}$ we see that $-\nabla E=-(\nabla C-\nabla(KC))=\dot KC+K\dot C+KCG=\dot KC$ on the one hand, and on the other hand $-\nabla E=-p^{N_a}(\nabla C-\nabla \tilde C)=\frac 1r\mathcal{E}$. As a consequence we have $\dot KC=\frac 1r\mathcal{E}$, or $K=\int \frac 1r\mathcal{E}C^{-1}d\G+\textrm{constant matrix}$. If we substitute $\G=0$, we see that this constant matrix will be the identity matrix. Hence we can conclude with
\[E=-\left(\int\frac 1r\,\mathcal{E}C^{-1}d\G\right)C.\]
Recall that $\eta=2g\log_p(g)+g$, then the remark after Proposition \ref{thm:boundsOnC} implies for $C^{-1}$ the bound
\[\ord\left((C^{-1})_k\right)\geq-\eta\lceil\log_p(k+1)\rceil.\]
The power series in the matrix $\frac 1r\mathcal{E}$ are integral, and integrating the coefficient of $\G^k$ adds at most an order of $\log_pk$. Finally we multiply with $C$, so that $\ord(E_k)\geq-(2\eta+1)\lceil\log_p(k+1)\rceil$.\ep

\noindent Remark that in fact we can prove more: let $D$ be a matrix such that $D\bmod\G=0$ and $\nabla D$ is integral. Then $D$ will satisfy the bounds for $E$ in the lemma.
\section{Complexity analysis}\label{sec:complexity}
All the time estimates in this section are measured in bit operations, and memory requirements are expressed in bit space.

\subsection{$p$-Adic computations}\label{sec:pAdicCalc}
In this section we work modulo $p^\nu$ for some positive integer $\nu$, and we assume that every $p$-adic number has valuation at least $-\O(\nu)$. Let $\Qq=\Qp[x]/\chi(x)$ as in step \ref{step:1}, then the memory requirements for an element of $\Qq$ are $a\nu\log_2p$, or $\O(a\nu)$ if we ignore the dependency on $p$. As pointed out in \cite{BernsteinFastMultiplication}, computing in such a quotient ring is essentially linear in the element size, hence $\Ot(a\nu)$.

Let us first consider calculating a power of Frobenius $\sigma^k$ with $k=\O(a)$ of an element $\alpha(x)$ of $\Qq$; we can suppose $\alpha(x)\in\Zq$ as $\sigma$ acts trivially on $\Qp$. We start with computing $\sigma^k(x)$, and in the next step $\sigma^k(\alpha(x))=\alpha(\sigma^k(x))$. We have that $\sigma^k(x)\equiv x^{p^k}\bmod p$ and $\chi(\sigma^k(x))=0$. First compute $x^{p^k}\bmod p$ by repeated squaring (this requires $\O(k)$ arithmetic operations in $\Fq$). It follows now from the irreducibility of $\bar{\chi}$ that $\chi'(x^{p^k})\not\equiv 0\bmod p$, so we can start quadratic Newton approximation. The total cost for $\sigma^k(\alpha(x))$ is $\Ot(a^2\nu)$, if we evaluate $\alpha$ at $\sigma^k(x)$ in a (naive) efficient way as with Horner's method, using $\O(a)$ operations in $\Zq\bmod p^\nu$. The memory cost is just $\O(a\nu)$.\\

In step \ref{step:7} we need a representation of $\Q_{q^n}$ as $\Qq[y]/\varphi(y)$ such that $y=\g$, a Teich\-m\"{u}l\-ler lift of $\og$. As described by Vercauteren \cite[Section 12.1.2]{CohenFrey} we can compute in time $\Ot(an\nu)$ a polynomial $f(z)$ --- called a Teichm\"{u}ller modulus --- over $\Qp$ such that $\Q_{q^n}=\Qp[z]/f(z)$ and $z=\g$. Here $\bar{f}$, the reduction of $f$ modulo $p$, is the minimal polynomial of $\og$ over $\Fp$ (to be computed in time $\Ot(a^2n^2)$ as in \cite{ShoupMinimalPolynomial}). Now we know that as $\varphi(z)=0$, we can split $f$ over $\Qq$ as $\varphi\cdot\varphi'$ for some polynomial $\varphi'$ such that $\bar{\varphi}$ and $\bar{\varphi'}$ are relatively prime. As it is easy to compute $\bar{\varphi}$ and $\bar{\varphi'}$, again using \cite{ShoupMinimalPolynomial}, we can use classical Hensel lifting \cite[Section 15.4]{ModernCompAlg} for computing $\varphi$ in time $\Ot(an\nu)$.\\

Finally we have to consider the action of Frobenius $\sigma^k$ on $\Q_{q^n}$. The dependency on $n$ of this step will dominate by far the overall complexity, hence any improvement here will immediately give better algorithms. Assuming $k=\O(an)$, we start again with computing $\sigma^k(y)$. A naive approach would be to reduce $y^{p^k}$ modulo $\varphi(y)$, but we can do this in a better way. Namely, compute in time $\Ot(a^2n^2)$ the projection $\bar{y}^{p^k}$ in $\F_{q^n}=\Fq[y]/\bar{\varphi}(y)$. Now we know that $\sigma^k(y)$ is also a Teichm\"{u}ller lift, and Section 12.8.1 of \cite{CohenFrey} shows how to compute $\sigma^k(y)$ in time $\Ot(an\nu)$. For an element $\alpha(y)\in\Q_{q^n}$, we can compute $\sigma^k(\alpha)$ in two different ways, resulting in the two complexity estimates of Theorem \ref{thm:firstResult}. We may assume that $\alpha(y)\in\Z_{q^n}$. Compute $\alpha^{\sigma^k}(y)$ in time $\Ot(na^2\nu)$, which is the degree $n-1$ polynomial over $\Zq$ obtained by applying $\sigma^{k\bmod a}$ to the coefficients of $\alpha$. Using Horner's method we obtain $\alpha^{\sigma^k}(\sigma^k(y))$ modulo $\varphi(y)$ in time $\Ot(n(an\nu))$ with memory requirements $\O(an\nu)$. All together we can thus compute $\sigma^k$ on an element of $\Q_{q^n}$ in time $\Ot(n^2a^2\nu)$ and space $\O(an\nu)$. A second method is more complicated but gives the faster algorithm of Theorem \ref{thm:firstResult}. We will treat it in the following section, where (\ref{eq:FastestEstimate}) gives the fastest result, and (\ref{eq:tradeOff}) gives for each $b\in[0,0.5]$ a trade-off between required memory and time.

\subsection{Fast modular composition of polynomials}\label{sec:modularComposition}
The following exposition is based on Section 12.2 of \cite{ModernCompAlg}, but we have added a certain trade-off in time versus memory.

Write $g(y):=\alpha^{\sigma^k}(y)$ and $\eta(y) := \sigma^k(y)$, then these are polynomials over $\Zq$ of degree at most $n-1$, whereas $\varphi(y)$ has degree $n$. Our goal is to compute $g(\eta(y))\bmod\varphi(y)$.

Let $b\in(0,0.5]$ and define $m := \lceil n^b\rceil$ and $m':=\lceil n/m\rceil\approx n^{1-b}$. We proceed in a number of steps.
\setcounter{step1}{0}
\begin{step1} We compute (trivially) polynomials $g_0(y),\ldots,g_{(m'-1)}(y)$ over $\Zq$ of degree at most $m-1$ such that
\[g(y)=g_0(y)+g_1(y)y^m+\cdots+g_{m'-1}(y)y^{(m'-1)m}\]
and compute $\eta(y)^i\bmod\varphi(y)$ for $i=0,\ldots,m-1$. This can be achieved in time $\Ot(man\nu)$ and space $\O(man\nu)$.
\end{step1}
\begin{step1}
For a polynomial $f(y)$ we denote by $[f(y)]$ the row matrix constructed from its coefficients (if necessary padded with zeroes), e.g.\ $[y+2y^3] = (0\ \ 1\ \ 0\ \ 2)$. Consider the following product of an $m'\times m$ with an $m\times n$ matrix over $\Zq$:
\begin{equation}\label{eq:modularCompMatrix}
\begin{pmatrix} [r_0(y)]\\ [r_1(y)] \\ \vdots \\ [r_{m'-1}(y)] \end{pmatrix} := \begin{pmatrix} [g_0(y)]\\ [g_1(y)] \\ \vdots \\ [g_{m'-1}(y)] \end{pmatrix}\cdot \begin{pmatrix} [\eta(y)^0]\\ [\eta(y)^1] \\ \vdots \\ \qquad[\eta(y)^{m-1}]\qquad \end{pmatrix}.
\end{equation}
Continue now with step 3 \emph{or} step 4.
\end{step1}
\begin{step1}
Taking $b=0.5$, $m\approx m'\approx \sqrt{n}$ we have to compute in (\ref{eq:modularCompMatrix}) the product of an $m\times m$ and an $m\times m^2$ matrix, and as proven in \cite[Section 5]{HuangPan} this can be done using $\O(m^{3.334})$ operations in $\Zq$. This yields a time complexity of $\Ot(n^{1.667}a\nu)$, and the memory requirements for (\ref{eq:modularCompMatrix}) are $\O(mna\nu)=\O(n^{1.5}a\nu)$. We can conclude that it is possible to compute $\sigma^k$ on an element of $\Z_{q^n}$ modulo $p^\nu$ in \begin{equation}\label{eq:FastestEstimate}\Ot(n^{1.667}a^2\nu) \text{ time and }\O(n^{1.5}a\nu)\text{ space.}\end{equation}
\end{step1}
\begin{step1}
If we do not compute the product (\ref{eq:modularCompMatrix}) at once, we can gain some memory. Note that representing the result of the product in (\ref{eq:modularCompMatrix}) requires already $\O(m'na\nu)$ bits of space. We divide the matrix formed by the polynomials $g_i(y)$ into $\lceil m'/m\rceil$ matrices of size $m\times m$: the first one, say $G_0$, is given by $\{g_0(y),\ldots,g_{m-1}(y)\}$, and generally $G_k$ is given by $\{g_{km}(y),\ldots,g_{km+m-1}(y)\}$. The corresponding matrix products are computed one by one, schematically
\[\begin{pmatrix} \begin{boxedminipage}[c]{2cm}\ \end{boxedminipage}\\\vdots\\\begin{boxedminipage}[c]{2cm}\ \end{boxedminipage}\,\end{pmatrix} = \begin{pmatrix} \begin{boxedminipage}[c]{\baselineskip}\ \end{boxedminipage}\\\vdots\\\begin{boxedminipage}[c]{\baselineskip}\ \end{boxedminipage}\,\end{pmatrix}\cdot \begin{pmatrix}\begin{boxedminipage}[c]{2cm}\ \end{boxedminipage}\, \end{pmatrix}.\]
Define $g^{(-1)} := 0$ and compute for every $k$
\begin{equation}\label{eq:inductiveg}
g^{(k)} := g^{(k-1)} + \sum_{i=0}^{m-1} r_{mk+i}(y)\cdot (\eta(y)^m)^i\mod \varphi(y).\end{equation}
This has to be alternated with the computation of the relevant matrix products above. At the end we have computed $g^{(\lceil m'/m\rceil -1)}=g(\eta(y))\bmod \varphi(y)$, but used only
\[\O(man\nu)=\O(n^{1+b}a\nu)\]
bits of memory, the requirement for each single matrix product. It is clear that (\ref{eq:inductiveg}) requires time $\Ot(man\nu)$ and has to be done approximately $m'/m$ times, resulting in $\Ot(m'a\nu)=\Ot(n^{2-b}a\nu)$. We will now investigate the time requirements for the matrix products. It is easy to verify that the product of an $m\times m$ with an $m\times m^c$ for some $c>1$ takes no more time than $\O(m^\omega (m^c/m))=\O(m^{c-1+\omega})$ operations in $\Zq$. We have to do this approximately $m'/m$ times, and as $c=1/b$ we find
\[\O(m^{\frac 1b -1+\omega}\frac {m'}m\,a\nu)=\O(n^{2+(3-\omega)b}a\nu)=\O(n^{2-0.624b}a\nu).\]
We note that for $b=0.5$ this gives $\O(n^{1.688}a\nu)$, precisely the same result as in \cite[Corollary 12.5]{ModernCompAlg}.

Combining these results with the $\Ot(a^2n^2)$ computation of $\bar y^{p^k}$ in Section \ref{sec:pAdicCalc}, we can compute $\alpha^{\sigma^k}(\sigma^k(y))$ in
\begin{equation}\label{eq:tradeOff}\text{time }\Ot(n^{2-0.624b}a^2\nu)\text{ and space }\O(n^{1+b}a\nu).\end{equation}
\end{step1}

\subsection{The complexity of the algorithm}\label{sec:complexitySteps}
We will use the following bounds: $N_0=\O(nga)$, $N_8$, $N_b$ and $M$ are all $\O(nga\log g)$, $N_{\G}=\O(ng^2a\kappa\log g)$, $N_a=\O(nga\kappa\log^3g)$ and the last constants $N_3$, $N_4$ and $N_6$ are $\Ot(g)$. For example, representing an element of $\Qq$ requires then $\O(aN_a)=\O(nga^2\log^3g)$ bits. We recall that we work modulo $\G^{N_{\G}}$, in steps 1 to 6 modulo $p^{N_a}$ and in steps 7 and 8 modulo $p^{N_b}$.
\setcounter{step1}{0}
\begin{step1}
We can compute $r$, $\alpha$ and $\beta$ with classical Gaussian elimination, which results in time $\Ot(g^3\cdot g\k\cdot aN_a)$.  Here we use an $\O(g\kappa)$ bound for the resulting degree in $\G$. The net result is $\Ot(ng^5a^2\kappa^2)$. We could gain one factor $g$ by using more sophisticated methods for computing the resultant as e.g.\ in \cite{MonaganResultants}.
\end{step1}
\begin{step1}
We have to apply formula (\ref{eq:ReductionXPowers}) $\O(g)$ times, and each time this requires $\Ot(g^2\kappa aN_a)$ calculations. Indeed, we can work with polynomials of $\deg_X\leq \O(g)$ and $\deg_\G\leq \O(g\kappa)$ as proven in Proposition \ref{thm:boundsOnH}. As we have $2g$ basis elements to consider, this gives together $\Ot(ng^5a^2\kappa^2)$.
\end{step1}
\begin{step1}
The computation of a single $C_k$ consists of at most $\O(g\kappa+\rho)$ matrix products and such a product requires $\Ot(g^\omega aN_a)$ time. Here $\omega$ is chosen bigger than the exponent for matrix multiplication as defined in \cite{ModernCompAlg}, and we can take $\omega=2.376$. The number of $C_k$'s to be computed is $N_{\G}$, so an overall time complexity is $\Ot(g^{1+\omega}a\k N_aN_{\G})=\Ot(n^2g^{4+\omega}a^3\kappa^3)$.

The memory requirements for representing $C$ dominate this step, and are $\O(aN_a N_{\G}g^2)=\O(n^2g^5a^3\kappa^2\log^3g)$. Note that for this last result we have to use Proposition \ref{thm:boundsOnC} to ensure that there appear no denominators that are too big. Steps 4 and 5 require certainly less memory space than this step.
\end{step1}
\begin{step1}
The calculation of $C^{\sigma}(\G^p)$ needs time $\Ot(g^2aN_aN_{\G})$, and the inversion costs $\Ot(g^\omega a N_aN_{\G})$.
\end{step1}
\begin{step1}
If we look at the analysis of Kedlaya and take into account the higher accuracy $N_a$, we find a time complexity of $\Ot(g^2N_a^2)=\Ot(g^4a^2n^2)$.
\end{step1}
\begin{step1}
Computing $F(\G)$ as a matrix of power series needs $\Ot(g^\omega aN_aN_{\G})$ time, and requires memory space of size $\O(g^2aN_aN_{\G})$. Multiplying $F(\G)$ and $r^M$ requires time $\Ot(g^2aN_aN_{\G})$.
\end{step1}
\begin{step1} This step is explained in more detail in Section \ref{sec:pAdicCalc}. To compute $\bar\varphi$ we can use Shoup's algorithm which requires $\Ot(a\sqrt{n}+n^2)$ time, and this can be ignored in the global complexity. To determine $\varphi$ we need $\Ot(anN_b)$ of time, and to reduce $g(y)$ modulo $\varphi(y)$ takes $\Ot(aN_bN_{\G})$ time due to \cite{BernsteinFastMultiplication}. Multiplying by $r(\g)^{-M}$ is negligible. As we have to do all this for the entire matrix, this step costs $\Ot(n^2g^5a^3\k)$ and the size of the resulting matrix --- with entries in $\Q_{q^n}$ --- is $\O(anN_bg^2)=\O(n^2g^3a^2\log g)$.
\end{step1}
\begin{step1}
Computing a $\O(an)$th power of Frobenius as explained in the previous section costs time $\Ot(n^{2.667}ga^3)$ respectively $\Ot(n^3ga^3)$ with memory requirements of $\O(n^{2.5}ga^2\log g)$ respectively $\O(n^2ga^2\log g)$. For determining the product of the matrices we use fast matrix multiplication, resulting in time $\Ot(g^\omega anN_b)=\Ot(n^2g^{1+\omega}a^2)$. Finally we compute the determinant using Newton's formula, for which we have to compute $\mathcal{F},\mathcal{F}^2,\ldots,\mathcal{F}^{2g}$ in time  $\Ot(g^{1+\omega}anN_b)=\Ot(n^2g^{2+\omega}a^2)$.
\end{step1}

Taking the maximum of all the above requirements we find immediately Theorem \ref{thm:firstResult}. For the second theorem we only need to keep $F(\G)$ modulo $p^{N_b}$ and to use steps \ref{step:7} and \ref{step:8}, so this theorem follows as well. The following picture gives a graphical overview of the time/memory requirements for our and Kedlaya's algorithm, where the full line indicates the trade-off explained in Section \ref{sec:modularComposition}. The values on the axes are to be seen as exponents for $n$ in the complexity.

\begin{center}
\setlength{\unitlength}{.3mm}
\begin{picture}(160,160)
\dottedline{2}(10,10)(35,10)\drawline(35,10)(150,10)
\drawline(35,10)(35,7)\drawline(85,10)(85,9)\drawline(135,10)(135,7)
\put(0,32){2}\put(0,132){3}\put(-7,155){\texttt{MEMORY}}
\put(10,140){\vector(0,1){10}}
\dottedline{3}(10,35)(150,35)\dottedline{3}(10,85)(150,85)\dottedline{3}(10,135)(150,135)

\drawline(10,150)(10,35)\dottedline{2}(10,35)(10,10)
\drawline(10,35)(7,35)\drawline(10,85)(9,85)\drawline(10,135)(7,135)
\put(33,-3){2}\put(133,-3){3}\put(155,-3){\texttt{TIME}}
\put(140,10){\vector(1,0){10}}
\dottedline{3}(35,10)(35,150)\dottedline{3}(85,10)(85,150)\dottedline{3}(135,10)(135,150)

\drawline(104,85)(135,35)\put(104,4){\tiny{2.688}}
\dottedline{5}(104,85)(104,10)
\put(100,85){\circle*{4}} \put(50,89){Deformation}
\dottedline{5}(100,85)(100,10)\put(98,-2){\tiny{2.667}}
\put(135,35){\circle*{4}} \put(139,39){Deformation}
\put(135,135){\circle*{4}} \put(139,139){Kedlaya}
\end{picture}
\end{center}
\section{Remarks}\label{sec:conclusion}
\subsection{The elliptic curve case}
For elliptic curves the above results are particularly interesting. It is well known that any elliptic curve in characteristic unequal to 2 can be described up to isomorphism by an equation
\[Y^2=X\cdot(X-1)\cdot(X-\lambda)\]
for some parameter $\lambda$; this is called the Legendre normal form. Moreover, if the original curve is defined over $\Fq$ and has $\Fq$-rational 2-torsion, then $\lambda$ will be in $\Fq$ as well. If we define $Q(X,\G):=X\cdot(X-1)\cdot(X-\G+1)$ over $\Fp$, we can hence reach every elliptic curve with rational 2-torsion with our deformation theory, and this over the prime field ($a=1$) and with a linear deformation ($\k=1$). Theorem \ref{thm:firstResult} gives then that we can compute the zeta function of an elliptic curve over $\F_{p^n}$ in time $\Ot(n^3)$ and space $\O(n^2)$ or time $\Ot(n^{2.667})$ and space $\O(n^{2.5})$. More details can be found in Section \ref{sec:elliptic}.

If we compare this with Harley's result \cite{HarleyECQuadratic}, we see that his algorithm has the same space complexity, but the time complexity is better, $\Ot(n^2)$. The main reason behind this difference is the following. We have to compute the `norm' of a $2\times 2$ matrix over $\Q_{p^n}$, which is the only step that requires time $\Ot(n^3)$. Harley however only has to compute the norm of an element of $\Q_{p^n}$, and for this he uses a classical formula that expresses this norm as a resultant. More precisely, let $\Z_{p^n}=\Zp[y]/f(y)$ with $f(y)$ a Teichm\"{u}ller modulus (see Section \ref{sec:pAdicCalc}), then for $\alpha(y)\in\Z_{p^n}$ we have
\[\mathcal{N}_{\Z_{p^n}/\Zp}(\alpha(y))=\text{Res}_y(\alpha(y),f(y)),\]
and this resultant can be evaluated very fast using an adaptation of  Moenck's \textsc{xgcd} algorithm. All details of Harley's result can be found in Section 3.10 of \cite{VercauterenThesis}. In \cite{HubrechtsQuadTime} we present a variant of our algorithm for genus 1 which has the same complexity as Harley's method. This works by `semi-diagonalising' the matrix $F(\g)$ and using the above fast norm computation. An implementation of this algorithm allows us to compute the zeta function of a random elliptic curve over $\F_{3^{100}}$ in less than one second.

If the curve has no $\Fq$-rational 2-torsion, the parameter $\lambda$ will be in an extension field of $\Fq$ of degree at most 3, and the curves are isomorphic over an extension field of degree at most 6 (see Proposition III.1.7 in \cite{Silverman}). If we compute the zeta function over this extension field, then it is possible to recover the original zeta function in an efficient but nondeterministic way, more details can be found in our forthcoming Ph.D. thesis \cite{HubrechtsThesis}. However, in our paper \cite{HubrechtsQuadTime} we were able to avoid this problem by choosing better families than the Legendre family.

\subsection{More than one parameter families?}
An interesting question would be to see how many hyperelliptic curves we can reach with a deformation as described in this paper. More precisely, for getting more curves over $\F_{q^n}$ we could increase the base field $\Fq$, the $\G$-degree $\k$ of $Q(X,\G)$, and even take the substitution $\G\leftarrow\og$ in a bigger field than $\F_{q^n}$. Although we have no decisive answer on this matter, the fact that the moduli space of hyperelliptic curves of genus $g$ in odd characteristic has dimension $2g-1$ suggests that this is not possible in a useful way for $g\geq 2$. Indeed, the number of `families of degree $\leq \k$' of the form $Y^2=Q(X,\G)$ is clearly at most $q^{(\k+1)(2g+1)}$, as there are $2g+1$ coefficients, all of the form $\alpha_{\k}\G^{\k}+\alpha_{\k-1}\G^{\k-1}+\cdots$ with all $\alpha_i\in\Fq$. Every family gives no more than $q^n$ curves, so an upper limit for the number of curves reached with a degree $\k$ deformation is $q^{(\k+1)(2g+1)}q^n$. The dimension of the moduli space gives that there exist about $q^{(2g-1)n}$ hyperelliptic curves of genus $g$ over $\F_{q^n}$, hence for $\k$ and $g\geq 2$ fixed  we really can have only a very small part of the possible curves. Of course it would be possible to increase $a$ and/or $\k$ to $\O(n)$, but this would not give a good algorithm (at all).
\section{Implementation results}\label{sec:immplementation}
Until now we have considered only theoretical complexities, but in order to be of any practical use, the algorithm has to work at a decent speed in a concrete implementation. We have implemented a version of the algorithm with theoretical time complexity $\Ot(n^3)$ and space $\O(n^2)$, using the computational algebra system Magma\footnote{See \texttt{http://magma.maths.usyd.edu.au/}.} V2.12-14.

Our implementation\footnote{Available on \texttt{http://wis.kuleuven.be/algebra/hubrechts/}.} does not use the basis $\{X^idX/\sqrt{Q}\}_i$ for $H_{MW}^-$, but instead we use $\{X^idX/\sqrt{Q}^3\}_i$, which gives --- as pointed out by Kedlaya \cite[Section 3.5]{KedlayaComputingZetaFunctions} --- an integral matrix of Frobenius. This gives an easier to implement and faster program (although not asymptotically). We require the base polynomial $Q(X,\G)$ to be defined over a prime field (of odd characteristic), and use a naive lift for $\Qq$ and $\Q_{q^n}$ and a Teichm\"{u}ller lift for $\g$. The situation $\G=0$ is handled by Harrison's implementation of Kedlaya's algorithm as it is built-in in Magma. We will present a few concrete results of this algorithm, all achieved on a Pentium IV running at 2.4 GHz, using 1.5 GB of memory with SuSE Linux 9.0 as operating system.

In trying to speed up the algorithm, we found that the crucial parameter is $M$, which depends on the convergence behavior of the entries of $F(\G)$. Experimentally a good value turned out to be $M = (3n/2+10)pg/3$, but some more tuning of this parameter may make the algorithm faster.

\subsection{Elliptic curves}\label{sec:elliptic}
As mentioned before, Harley's algorithm has complexity $\Ot(n^2)$, and will hence be faster than ours, but to our knowledge it has not yet been implemented in Magma for odd characteristic --- Magma uses the SEA-algorithm for such curves\footnote{An implementation in C++ of a faster algorithm is described in \cite{Madsen}.}. We always work with the Legendre family and a random parameter in the finite field. For small fields SEA is faster than the deformation, but for field sizes starting from about $3^{100}$, $5^{60}$ or $7^{40}$ the deformation algorithm is substantially faster. A few timing results are gathered in the following table. Remark that all times are in seconds.

\begin{center}{\small
\begin{tabular}{|r||r|r|r|r|r|r|}
\hline $p\backslash n$\vphantom{$\sum^j$} & 20 & 40 & 60 & 80 & 100 & 150\\\hline\hline
3\vphantom{$\sum^j$} & 0.95 & 5.29 & 10.74 & 27.75 & 40.28 & 178\\\hline
5\vphantom{$\sum^j$} & 2.40 & 8.01 & 16.37 & 39.18 & 62.79 & 245\\\hline
29\vphantom{$\sum^j$} & 22.59 & 69.98 & 145.35 & 279.90 & 419 & 1300\\\hline
107\vphantom{$\sum^j$} & 158.67 & 472.52 & 1025.39 & 1742.27 & 2676 & 7669\\\hline
233\vphantom{$\sum^j$} & 520.64 & 1619.10 & 3379.69 & & & \\\hline
\end{tabular}}
\end{center}

\noindent A parameter $\og\in\F_{3^{500}}$ used 4605 seconds, and to get an idea of the object sizes involved: the $p$-adic accuracy in the beginning of the algorithm was $N_a=298$, the power series in $\G$ were computed modulo $\G^{N_{\G}}$ with $N_{\G}=3048$, and $M$ was equal to $760$. Most of the time was consumed by computing $\g$, $F(\g)$ and $\mathcal{F}$. It is indeed quite interesting to study which steps in the algorithm take most of the time. In the following table $\bar{F}(\G)$ stands for the analytical continuation of $F(\G)$ (which comes down to the product $r^M\cdot F(\G)$), and $C_\sigma^{-1}$ is short for $C^\sigma(\G^p)^{-1}$.

\begin{center}{\small\begin{tabular}{|l||r|r|r|r|r|r|r||r|r|r|}
\hline \vphantom{$\sum^{j^j}$}& $r(\G)$ & $H$ & $C$ & {$C_\sigma^{-1}$} & $F(0)$ & $F(\G)$ & $\bar{F}(\G)$ & $\g$ & $F(\g)$ & $\mathcal{F}$\\\hline\hline
$\F_{3^{50}}$\vphantom{$\sum^j$} & 0.26 & 0 & 0.18 & 0.06 & 0.08 & 1.19 & 0.76 & 0.52 & 1.38 & 3.11 \\\hline
$\F_{3^{100}}$\vphantom{$\sum^j$} & 0.26 & 0 & 0.37 & 0.13 & 0.14 & 3.87 & 2.26 & 2.65 & 6.71 & 23.76 \\\hline
$\F_{3^{150}}$\vphantom{$\sum^j$} & 0.25 & 0 & 0.55 & 0.19 & 0.20 & 8.39 & 5.37 & 12.23 & 25.09 & 125.5 \\\hline
$\F_{3^{500}}$\vphantom{$\sum^j$} & 0.26 & 0 & 2.36 & 0.83 & 0.93 & 93.03 & 55.46 & 280.9 & 424.9 & 3710 \\\hline
$\F_{7^{100}}$\vphantom{$\sum^j$} & 0.30 & 0 & 0.90 & 0.15 & 0.34 & 15.66 & 10.67 & 6.22 & 19.76 & 31.61 \\\hline
\end{tabular}}\end{center}

This table learns us that if the extension degree increases, the last steps (being precisely those which depend on $\og$) in the algorithm become more time consuming. This is to be expected, as they are precisely the steps cubic in $n$. Checking the correctness of the result took at most a few seconds.

To get an idea of the memory requirements, Magma reported 11.87 MB for $\g\in\F_{3^{100}}$ and for $\g\in\F_{3^{500}}$ the amount was 48 MB. This includes the kernel memory, approximately 3 MB.

\subsection{Higher genus}
The dependency on the genus is not as good as in Kedlaya's algorithm, and for $g\geq 5$ the algorithm is not very useful anymore. In particular the use of memory --- although being $\O(n^2)$ --- grows then too big for even very small $n$. We tested a few higher genus curves as follows: $Q$ is defined to be equal to $X^{2g+1}$ plus (a random polynomial over $\Fp$ of degree at most $2g$) plus $\G$ $\times$ (such a polynomial). For the parameter $\og$ again a random element is chosen. For these higher genus situations it is interesting to take advantage of the deformation in order to compute more zeta functions within a family. Therefore we present all timing results as $x/y$ where $x$ is the time needed for the precomputation, and $y$ the time for one parameter. The second columns below give the memory use in MB.

\begin{center}{\small
\noindent\begin{tabular}{|r||r|r||r|r||r|r|}
\hline $p^n\backslash g$\vphantom{$\sum^j$} & \multicolumn{2}{|c||}{2} & \multicolumn{2}{|c||}{3} & \multicolumn{2}{|c|}{4}\\\hline\hline
$3^{50}$\vphantom{$\sum^{j^j}$} & 165/45 & 18.88 & 1840/219 & 61 & 13899/745 & 175 \\\hline
(K) $3^{50}$\vphantom{$\sum^{j^j}$} & 166 & 24.28 & 772 & 84 & 2070 & 142 \\\hline
$3^{100}$\vphantom{$\sum^{j^j}$} & 520/281 & 28.10 & 7062/1402 & 109 & 46672/4898 & 356 \\\hline
(K) $3^{100}$\vphantom{$\sum^{j^j}$} & 1216 & 110.10 & 6104 & 399 & 17483 & 723 \\\hline
$3^{200}$\vphantom{$\sum^{j^j}$} & 2050/2120 & 54.29 & 28144/12026 & 238 & & \\\hline
(K) $3^{200}$\vphantom{$\sum^{j^j}$} & 11785 & 635.34 & & & & \\\hline
$3^{400}$\vphantom{$\sum^{j^j}$} & 8635/20586 & 122.54 & & & & \\\hline\hline
$5^{100}$\vphantom{$\sum^{j^j}$} & 1177/490 & 39.73 & 17356/2529 & 188 & 99886/8826 & 578 \\\hline
(K) $5^{100}$\vphantom{$\sum^{j^j}$} & 2796 & 221.30 & 16082 & 606 & 46842 & 1280 \\\hline
$31^{100}$\vphantom{$\sum^{j^j}$} & 23989/3817 & 273.30 & &  &  &  \\\hline
\end{tabular}}\end{center}

The (K) in front of a row means that Kedlaya's algorithm was used for precisely the same curve. The most striking difference is the case with genus 2 over $\F_{3^{200}}$, where we tried a couple of different random equations, all of which gave similar results\footnote{If we choose the equation such that the degree of the resultant is lower than 8, e.g.\ 5 or 4, than we get timings like 650/1650 and a memory use of 35 MB.}. We note also that it was not even possible to try a genus 2 curve over $\F_{3^{400}}$ with Kedlaya's method, as it would require approximately 5 GB of memory, eight times the amount used for $3^{200}$. For higher genus the advantage of deformation is more on the level of memory requirements\footnote{For $g=4$ and field size $3^{100}$ we could even gain about 100 MB by computing the entries of $F(\g)$ one by one given $C$, $C^\sigma(\G^p)^{-1}$ and $F(0)$. Of course the advantage of the family is lost this way.} and computing within families.\\

The conclusion is that for low genus and a field size big enough, a deformation algorithm can give a substantial advantage over the `classical approach', at least for curves in certain one parameter families. In particular the memory requirements drop dramatically, and after some precomputations it is rather efficient to compute concrete zeta functions.
\fontsize{9}{11pt}\selectfont
\bibliographystyle{acm}
\bibliography{bibliography}

\end{document}